\documentclass[titlepage,11 pt]{article}

\usepackage{amsmath} %Never write a paper without using amsmath for its many new command
\usepackage{graphicx} %If you want to include postscript graphic
\usepackage[english]{babel} % per la lingua inglese
\usepackage{amsfonts} %% per usare mathds

\usepackage{amsthm}

\usepackage{setspace}
\usepackage{lipsum}
\usepackage[margin=3cm]{geometry}

\usepackage[affil-it]{authblk}

\usepackage[authoryear]{natbib}
\usepackage{url} % not crucial - just used below for the URL

\usepackage{hyperref}
\usepackage{bbm}

\newcommand{\mf}{\mathbf}
\newcommand{\vf}{\mathbf}

\usepackage[T1]{fontenc}
\usepackage{a4wide, url}
\usepackage[english]{babel}
\usepackage{graphicx}
\usepackage{array}
\usepackage{multirow}
\usepackage{amsthm}
\usepackage{amsmath}
\usepackage{amsfonts}
\usepackage{amssymb}
\usepackage{pgf}
\usepackage{bm}
\usepackage{paralist}
\usepackage{pdflscape}
\usepackage{xfrac}
\newcolumntype{P}[1]{>{\centering\arraybackslash}p{#1}}
\usepackage{bbm}

\usepackage{chngcntr}
\usepackage{apptools}

\AtAppendix{\counterwithin{Ass}{section}}
\AtAppendix{\counterwithin{Lemma}{section}}
\AtAppendix{\counterwithin{Thm}{section}}
\AtAppendix{\counterwithin{Coroll}{section}}
\AtAppendix{\counterwithin{Def}{section}}

\newcommand{\numberset}{\mathbb}

\newcommand{\R}{\numberset{R}}

\newtheorem{Thm}{Theorem}%[section]
\newtheorem{Coroll}{Corollary}%[section]
\newtheorem{Lemma}{Lemma}%[section]
\newtheorem{Ass}{Assumption}%[section]
\title{\textbf{Large factor model estimation \\  by nuclear norm plus $l_1$ norm penalization}}

\author{\textbf{Matteo Farn\'{e}}  \thanks{Electronic address: \texttt{matteo.farne@unibo.it}; Corresponding author}}
\affil{Department of Statistical Sciences,\\ University of Bologna, Italy}

\author{\textbf{Angela Montanari}}  %&\thanks{Electronic address: \texttt{matteo.farne2@unibo.it}; Corresponding author}}
\affil{Department of Statistical Sciences,\\ University of Bologna, Italy}

\begin{document}

\maketitle

\begin{abstract}
This paper provides a comprehensive estimation framework via nuclear norm plus $l_1$ norm penalization for high-dimensional approximate factor models with a sparse residual covariance. The underlying assumptions allow for non-pervasive latent eigenvalues and a prominent residual covariance pattern. In that context, existing approaches based on principal components may lead to misestimate the latent rank, due to the numerical instability of sample eigenvalues. On the contrary, the proposed optimization problem retrieves the latent covariance structure and exactly recovers the latent rank and the residual sparsity pattern. Conditioning on them, the asymptotic rates of the subsequent ordinary least squares estimates of loadings and factor scores are provided, the recovered latent eigenvalues are shown to be maximally concentrated and the estimates of factor scores via Bartlett's and Thompson's methods are proved to be the most precise given the data. The validity of outlined results is highlighted in an exhaustive simulation study and in a real financial data example.
\end{abstract}

%%INTRODUCTION

\section{Introduction}

The digital revolution has enormously enlarged
the amount of available data for researchers and practitioners.
Consequently, the need rises to develop methodologies able to summarize the content of high-dimensional datasets, in order to derive meaningful information from them.

The factor model is an effective tool to this end,
as it detects the latent covariance structure behind a set of variables.
We can define the factor model for any $p$-dimensional mean-centered
random vector $x$ as %the following representation:
\begin{equation}x= B f+\epsilon,\label{mod_base}\end{equation}
where $B$ is a $p\times r$ matrix,
%such that $B'B=\Lambda_L$,
%$\Lambda_L$ is a $r\times r$ diagonal matrix,
$f$ is a $r \times 1$ random vector with
$E(f)=0_r$ and $V(f)=I_r$, and %$ X, {\epsilon}$
$\epsilon$ is a $p \times 1$ random vector with $E(\epsilon)=0_p$ and $V(\epsilon)=S^{*}$,
with $S^{*}$ full rank $p \times p$ matrix.
%It follows that $E(x)=0_p$, without loss of generality.
%Assuming that $f$ and $\epsilon$ are componentwise uncorrelated,
%we obtain $V(x)=\Sigma^{*}=L^{*}+S^{*}$, where $L^{*}=BB'$.
%We also set $B=U_L \Lambda_L^{1/2}$, where $U_L$ is a $p\times r$ matrix containing as columns the eigenvectors of $L^{*}$ and $\Lambda_L$ is the $r\times r$ diagonal matrix of the corresponding eigenvalues.
%, with $B=U_L \Lambda_L^{1/2}$.

Let us indicate by $\Sigma^{*}$ the $p \times p$ covariance matrix of the random vector $x$.
Assuming that $f$ and $\epsilon$ are componentwise uncorrelated,
the factor model (\ref{mod_base})
induces in $\Sigma^{*}$ a low rank plus residual decomposition of the following type: \begin{equation}\Sigma^{*}=L^{*}+S^{*}=BB'+S^{*},\label{mod_l+s}\end{equation}
%where $B$ is any $p\times r$ matrix such that
where $L^{*}=BB'=U_L \Lambda_L U_L'$,
with $U_L$ $p\times r$ semi-orthogonal matrix %containing
%as columns the eigenvectors of $L^{*}$
and $\Lambda_L$ $r\times r$ diagonal matrix. %of the corresponding eigenvalues.
%and $S^{*}$ is a $p \times p$ matrix with full rank. % semi-orthogonal %diagonal
The representation (\ref{mod_l+s}) is invariant under orthogonal transforms,
and it is therefore unidentifiable from the data without further constraints.

Suppose that we have a sample $x_k$, $k=1,\ldots,n$.
The unbiased $p \times p$ sample covariance matrix
is defined as $\Sigma_n=(n-1)^{-1}{\sum_{k=1}^n x_k x_k'}$.
Most of factor model estimation methods
rely on $\Sigma_n$ as an input,
and make essentially use of two techniques:
principal component analysis (see \cite{jolliffe2002principal} for an overview)
and maximum likelihood.
As outlined in \citep{bai2008large}, however, a large dimension $p$
leads to some particular estimation problems for model (\ref{mod_base}),
due to the limitations of $\Sigma_n$ in high dimensions.
%especially about the relationship between data dimensionality $p$ and sample size $n$.
%from the factor model perspective.

From a historical perspective, the classical inferential theory for factor models \citep{anderson1958introduction}
prescribes that the dimension $p$ is fixed while the sample size $n$ tends to infinity.
In particular, the strict condition $p<n$
is required to ensure consistency.
%the positive definiteness of the sample covariance matrix $\Sigma_n$.
As a consequence, the classical framework is clearly inappropriate if $p$ is large.
When $p>n$, in fact, $\Sigma_n$ becomes
inconsistent and no longer Wishart-distributed.
%thus making any estimate based on $\Sigma_n$ unusable.
%%1958 libro!

%%modello fattoriale subito!
%%

At the same time, when the dimension $p$ and the sample size $n$ are finite,
\cite{anderson1956statistical} show that
the use of the principal components of $\Sigma_n$ to estimate $B$ leads to
loadings and factor scores estimates which are incoherent with model assumptions,
because any estimate of $S^{*}$ so derived will never be full rank.
That is the reason why \cite{chamberlain1983arbitrage} prove that the principal components of $\Sigma_n$ consistently identify $L^{*}$ under model (\ref{mod_base}) as $p\rightarrow \infty$, provided that the $r$ eigenvalues of $L^{*}$ diverge with $p$ and $S^{*}$ is a non-diagonal matrix with vanishing eigenvalues as $p$ diverges.
%The residual $S$ is allowed to be non-diagonal and is estimated as a consequence.

Another relevant aspect concerns the ratio $p/n$.
If $p/n \rightarrow 1^{-}$,
the bad conditioning properties of $\Sigma_n$ inevitably affect the consistency
of principal component analysis (PCA) as a factor model estimation method.
In fact,
the sample eigenvalues follow the Marcenko-Pastur law \citep{marvcenko1967distribution},
which crucially depends on the ratio $p/n$.
In particular, if $p/n \rightarrow 1^{-}$, it is more likely to observe small sample eigenvalues,
thus making $\Sigma_n$ numerically unstable.

%enough to be caught.
%% No scree plot!

%via the standard criteria based on sample eigenvalues like those of baing02, because the smallest eigenvalues are likely not to be large enough.

An overall inferential theory of PCA
as a high-dimensional factor model estimation method has been developed in \cite{bai2003inferential}.
As also outlined in \cite{chamberlain1983arbitrage},
\cite{bai2003inferential} shows that the pervasiveness of the eigenvalues of $L^{*}$
as $p \rightarrow \infty$ is crucial for the exact recovery
of the latent rank $r$, performed by the identification criteria of \cite{bai2002determining}.
If that condition is violated, the latent rank $r$ may be underestimated by any PCA-based method,
as one or more latent eigenvalues may be unrecovered,
because the corresponding sample eigenvalues may not be large enough.
In order to achieve consistency,
PCA tolerates a non-diagonal residual covariance matrix $S^{*}$
and residual heteroscedasticity,
provided that $p$ and $n$ are both large and ${\sqrt{n}}/{p}$ tends to $0$.
On the contrary, if only $n$ is large, no non-diagonal residual covariance structure is admitted.

%% parte da covarianza!
\cite{fan2013large} propose to estimate the covariance matrix $\Sigma^{*}$ in high dimensions under representation (\ref{mod_l+s}) by taking out the principal components of $\Sigma_n$ and then thresholding their orthogonal complement, under the assumption that $S^{*}$ has a bounded $l_1$ norm as $p$ diverges. The uniform parametric consistency of loadings, factor scores and common components obtained by such covariance matrix estimates is established.
% assuming the residual covariance matrix $S^{*}$ to be sparse
%This procedure allows to relevantly decrease the number of parameters by imposing that sparsity assumption to .
That sparsity assumption on $S^{*}$ also allows to make the estimation error of $\Sigma_n$ vanish in relative terms as $p$ diverges. %, thus turning the curse to the blessing of dimensionality.
%In particular, it is required that the vector-induced $l_1$ norm of $S^{*}$ vanishes as $p$ diverges, while the eigenvalues of $L^{*}$ diverge with $p$.

The asymptotic distribution of factors and factor loadings estimated via PCA when both $p$ and $n$ are large is derived in \cite{bai2013principal}. A relevant merit of that paper is that factors and loadings are precisely identified without the need of any rotation.
Under relatively weak factors in terms of explained variance proportion, \cite{onatski2011asymptotic}
derives the (normal) asymptotic distribution of the coefficients in the OLS regressions of the PC estimates of factors (loadings) on the true factors (true loadings). That distribution has good approximation properties even when both $p$ and $n$ are reasonably small.

%A method to derive by principal components asymptotically consistent forecasts from many predictors is derived in \cite{stock2002forecasting} and applied in \cite{stock2002macroeconomic}. Principal components are also effectively used to estimate the generalized dynamic factor model in \cite{forni2000generalized}.

%blessing of dimensionality!

Concerning maximum likelihood estimation,
\cite{anderson1958introduction} shows that the exact maximum likelihood is consistent for loading estimation,
even if it is still inconsistent as far as factor scores estimation is concerned.
Nevertheless, factor scores can be consistently estimated by the conditional maximum likelihood,
via a frequentist approach (Bartlett's estimator) or a Bayesian approach (Thompson's estimator).

The consistency of maximum likelihood (ML) to estimate a high-dimensional factor model
has been studied in \cite{bai2012statistical} (previous contributions on the topic also include \cite{joreskog1967some} and \cite{lawley1971factor}).
Differently from the estimator of factor scores based on PCA, the one based on ML is consistent also for small $p$ and $n$,
even if the estimator distribution is less complicated to derive when $p$ diverges.
%, for the blessing of dimensionality.
ML has a better asymptotic rate and is more efficient than PCA in the case of independent and heteroscedastic residuals. However, in presence of a non-diagonal residual covariance structure, the convergence rates and the optimality conditions of ML estimators become cumbersome. It is important to note that the relative magnitude of $p$ and $n$
is a crucial issue for both methods (ML and PCA) to provide consistent factor model estimates.

%The inferential theory for large dimensional factor models was developed in Bai (2003), where the estimation is via PCA, and Bai and Li (2012), where the estimation is via ML. The PCA estimator tolerates serial residual correlation and heteroscedasticity only if $p$ and $n$ are both large and $\frac{\sqrt{n}}{p}$ tends to $0$. Otherwise, if $n$ is large, it requires a residual identity covariance. In any case, for small $p$ PCA is ineffective. On the contrary, the ML estimator works also when $p$ and $n$ are small, has a better rate and it is more efficient than PCA in case of independent and heteroscedastic residual, even if its distribution is much easier if $p$ is large (blessing of dimensionality). The case of serial correlation complicates rates and optimality conditions. For both estimators, asymptotic normality is ensured via appropriate CLTs, while the overall number of parameters represents an issue from theoretical and empirical point of view.

Given these premises, the interest arises to find an alternative estimation method to ML and PCA,
as they both present some relevant drawbacks in high dimensions. First of all,
the latent rank recovery fails if the latent eigenvalues are not spiked enough with respect to the dimension.
Then, the sample covariance matrix is increasingly numerically unstable as the dimension increases,
such that the need to regularize sample eigenvalues rises. In addition,
a more effective sampling theory is needed with respect to the degree of spikiness of
latent eigenvalues and the degree and pattern of residual sparsity.
Ideally, all these features should be present
also for finite values of $p$ and $n$.

%%Then,
%%inverte!

%In our days, some relevant issues are still open. First, rank selection: in case of non-pervasive eigenvalues, which is very frequent in nature, the latent rank is lost by baing02 criteria. For this reason, we replace PCA by the nuclear norm heuristics, that finds out automatically the latent rank in a consistent way, only provided that the smallest nonzero eigenvalue is large enough. The same holds for the minimum nonzero absolute entry of the sparse component. Differently from baing17, we discuss the algebraic consistency of the estimates, their optimality properties under certain conditions, and we show the systematic gain of our estimates with respect to POET ones.

In \cite{bai2019rank}, it is proposed to use the nuclear norm heuristics in place of PCA. That work provides the asymptotic normality and parametric consistency of approximate factor model estimates as both $p$ and $n$ diverge.
The proposed objective function is a least squares loss penalized by a nuclear norm plus $l_1$ norm heuristics,
which is useful to detect covariance matrix decompositions of type (\ref{mod_l+s})
where $S^{*}$ is element-wise sparse. In \cite{farne2020large},
the authors exploit the same heuristics to derive algebraically consistent covariance matrix estimates,
that is, the latent rank and the residual sparsity pattern
are exactly recovered for finite values of $p$ and $n$.
%under their framework.
Such a feature is extremely important, as it allows to avoid the use of any identification criterion for the latent rank like the one described in \cite{bai2019rank}.

The results of \cite{farne2020large} are obtained by allowing for intermediate degrees of spikiness
for latent eigenvalues and intermediate degrees of sparsity for the residual component.
In particular, their assumptions prescribe
that the latent eigenvalues are spiked
in the sense of Yu and Samworth (\cite{fan2013large}, p. 656),
thus allowing for intermediately pervasive latent factors as $p$ diverges.
%\begin{Ass}\label{eigenvalues}
%All the eigenvalues of the $r\times r$ matrix $p^{-\alpha}\mf{B}^\top \mf{B}$
%are bounded away from $0$ for all $p$ and $\alpha \in [0,1]$.
%\end{Ass}
What is more, the number of non-zeros in the residual component $\mf{S}^{*}$ is allowed
to grow with $p$ (even if slower than the latent eigenvalues).
The identifiability of
the matrix components $\mf{L}^{*}$ and $\mf{S}^{*}$
is ensured by imposing that $\mf{L}^{*}$ and $\mf{S}^{*}$ are far enough
from being sparse and low rank respectively.
We refer to Appendix \ref{ass:fm20} for technical details.

In this paper, we provide the finite error bounds
for loadings, factor scores and common components estimated under the framework of \cite{farne2020large}.
The theoretical background %, \far{which builds upon Assumptions \ref{eigenvalues} and \ref{alg}},
is discussed in Section \ref{sec:back}. %under which both POET and UNALCE estimates are consistent.
In Section \ref{sec:ols_both}, the asymptotic consistency of factor model estimates based on
the nuclear norm plus $l_1$ norm heuristics under those conditions is proved.
In Section \ref{sec:unalce},
we present a re-optimized version of those estimates,
from which in Section \ref{sec:eig_conc} the most precise factor model estimates
produced by any algebraically consistent low rank and sparse component estimates
%from a numerical point of view
are derived, given the data.
In Section \ref{sec:bt_fac}, we highlight that the subsequent Bartlett's and Thompson's estimators of factor scores provide the tightest error bound in Euclidean norm
within the classes of algebraically consistent low rank and sparse component estimates, given the data.
In Section \ref{sec:sim}, we provide a wide simulation study proving the validity of our approach.
Section \ref{sec:real} then shows a real financial data application.
Finally, the conclusions follow in Section \ref{sec:concl}.
%(Bartlett's and Thompson's)

%most stable
%Euclidean norm

%The paper proceeds as follows. ...
%In Section 2, we establish asymptotic rates for our estimates following the framework of \cite{bai2003inferential} under the assumptions of \cite{farne2018finite}. In Section 3, we prove that our estimates have the most possible concentrated eigenvalues into the classes of algebraically consistent solutions. In Section 4 we describe a simulation study in support of our results. Section 5 reports conclusions and discussions.

%%BACKGROUND
%%L+S principal components.

\section{Theoretical background}\label{sec:back}

\subsection{Notation}

%Before stating them, we need to introduce some matrix norm notation.
Given a $p\times p$ symmetric positive-definite matrix $M$,
we denote by $\lambda_i(M)$, $i\in\{1,\ldots, p\}$
the eigenvalues of $M$ in descending order.
Then, we recall the following norm definitions:
\begin{description}
\item{(i)} Element-wise:
\begin{description}
\item{(a)} $L_0$ norm: $\Vert M\Vert_0=\sum_{i=1}^p \sum_{j=1}^p \mathbbm{1}( M_{ij}\ne 0)$, which is the total number of non-zeros;
\item{(b)} $L_1$ norm: $\Vert M\Vert_1=\sum_{i=1}^p \sum_{j=1}^p \vert M_{ij}\vert$;
\item{(c)} Frobenius norm: $\Vert M\Vert_{F}=\sqrt{\sum_{i=1}^p \sum_{j=1}^p  M_{ij}^2}$;
\item{(d)} Maximum norm: $\Vert M\Vert_{\infty}=\max_{i\leq p, j \leq p} \vert M_{ij}\vert$.
\end{description}
\item{(ii)} Induced by vector:
\begin{description}
\item{(a)} $\Vert M\Vert_{0,v}=\max_{i \leq p} \sum_{j \leq p} \mathbbm{1}( M_{ij} \ne 0)$,
which is the maximum number of non-zeros per column,
defined as the maximum `degree' of $ M$;
\item{(b)} $\Vert M\Vert_{1,v}=\max_{i \leq p} \sum_{j \leq p} \vert M_{ij}\vert$;
\item{(c)} Spectral norm: $\Vert M\Vert_{2}=\lambda_1( M)$.
\end{description}
\item{(iii)} Schatten:
\begin{description}
\item{(a)} Nuclear norm of $ M$, here defined as the sum of the eigenvalues of $ M$:
$\Vert M\Vert_{*}=\sum_{i=1}^p \lambda_i(M)$.%, which is.
\end{description}
\end{description}
Given a $p-$dimensional vector $v$, we denote by $\Vert v \Vert=\sqrt{\sum_{i=1}^p v_i^2}$ the Euclidean vector norm of $v$.

%`modello all'inizio!
% radice di \Lambda_L!

%%logica
%%linguaggi

%$f_i$ and $\epsilon_i$ are $r \times 1$ and $p \times 1$ zero-mean random variables.
%with finite fourth moments.
%The covariance matrix of $f_i$ is $I_r$, the one of $\epsilon_i$ is $S^{*}$.  $S^{*}$
%is assumed to have only $s \ll p(p-1)/2$ non-zero elements.

\subsection{State of the art}

%We observe a sample of $p$-dimensional data vectors $x_i$, $i=1,\ldots,n$.
Imposing $S^{*}=I_p$, \cite{bai2003inferential} shows that the loading matrix $B$ and the factor scores $f_k$, $k=1,\ldots,n$, are consistently recovered under model (\ref{mod_base}) as $p\rightarrow \infty$
by extracting the top $r$ eigenvectors of $\Sigma_n$,
provided that the $r$ eigenvalues of $L^{*}$
%, differently from the $p$ eigenvalues of $S^{*}$,
are scaled with $p$.
%(see Section \ref{sec:ols_both}).
%are consequently recovered at each $k$ via conditional least squares.
The reason why this method is consistent as $p\rightarrow \infty$ can be understood recalling \cite{hotelling1933analysis}. In fact, the principal components of $\Sigma_n$ are derived solving the problem \begin{equation}\min_{L, rk(L)\leq r} ||\Sigma_n-L||_{F},\end{equation}
which is equivalent to the problem \begin{equation}\min_{B_{j},f_{k,j}} \frac{1}{n} \sum_{k=1}^n ||x_k-z_k||_2,\label{pca2}\end{equation}
where $z_k=\sum_{j=1}^r B_{j} f_{k,j}$, $f_{k,j}$ is the $j$-th component of $f_k$
and the column vectors $B_j$, $j=1,\ldots,r$, are orthogonal. %at different $j$.
Intuitively, the solutions to problem (\ref{pca2}) are consistent under model (\ref{mod_base})
if and only if the eigenvalues of $L^{*}$ are scaled with $p$ and $p\rightarrow \infty$,
because otherwise the signal $z_k$ would not be strong enough to be detected.

Full solution vectors $B_j$, $j=1,\ldots,r$, can be difficult to interpret in high dimensions.
For this reason, \cite{zou2006sparse} introduce the sparse principal component analysis (SPCA), a method based on a version of problem (\ref{pca2}) where each $B_j$ is penalized by a ridge plus lasso penalty. The resulting sparse principal components are not orthogonal anymore and represent approximate solutions,
which reduce effectively the complexity of estimated components when $p$ is large.

At the same time, as $p$ diverges, the assumption $S^{*}=I_p$ is definitely too strong, as it is unlikely that the latent structure is able to entirely catch the covariance for all pairs of variables.
In order to relax that assumption, \cite{candes2011robust}
propose the principal component pursuit (PCP),
that is based on the solution of the following problem:
\begin{equation}\min_{L+S=\Sigma_n} ||L||_{*}+ ||S||_1,\label{pcp}\end{equation}
where $||L||_{*}$ is the nuclear norm of $L$, which is the sum of its singular values, and $||S||_1$ is the $l_1$ norm of $S$, which is the sum of all its absolute entries. Problem (\ref{pcp}) can be thought of as a robust PCA problem in presence of missing or grossly corrupted data.
It is solved exploiting the singular value thresholding algorithm of \cite{cai2010singular}.

Even if problem (\ref{pcp}) is able to bypass the assumption $S^{*}=I_p$,
the number of parameters to be recovered may be remarkably high
without any further assumption on $S^{*}$, particularly if $p$ is large.
In order to reduce the parameter space dimensionality, a rough alternative is to impose sparsity on $\Sigma^{*}$. In the covariance matrix context, for instance, \cite{bickel2008covariance} assume that $\Sigma^{*}$ is sparse and recover it by solving the problem $\min_{\Sigma} ||\Sigma_n-\Sigma||_1$. This problem is solved by applying to $\Sigma_n$
the soft-thresholding algorithm of \cite{daubechies2004iterative}, which
is consistent for $\Sigma^{*}$ but does not provide any dimension reduction.

%%zou et al.

% modello
% metodi di stima
% inconveienti
% nuovo!

%A successful attempt to merge the two aspects (sparsity and dimension reduction) was done

The use of the nuclear norm for rank minimization as an alternative to PCA was first proposed in \cite{fazel2001rank}.
The nuclear norm was then successfully applied to matrix completion problems, among which the Netflix problem is the most celebrated one. Within this research strand, we mention \cite{srebro2005maximum}, \cite{candes2010power}, \cite{mazumder2010spectral}, and \cite{hastie2015matrix}, which all describe and solve approximate robust PCA problems. %In addition, the nuclear norm has been used to estimate the covariance matrix under the low rank plus sparse assumption in \cite{agarwal2012noisy}.

Given these premises, in this paper we merge dimension reduction and sparsity in a single problem with the aim to explore the performance of the subsequent estimates of factor scores and loadings.
First, we recover the two components $L^{*}$ and $S^{*}$ of $\Sigma^{*}$ from $\Sigma_n$.
This step is performed by solving the following problem \citep{farne2020large}:
\begin{equation}
\min_{L,S} ||\Sigma_n-(L+S)||_{F}+\lambda ||L||_{*}+\rho ||S||_{1,off},\label{ours}%
\end{equation}
where $||L||_{*}$ is the nuclear norm of $L$ and $||S||_{1,off}$ is the $l_1$ norm of
$S$ excluding the diagonal, i.e. $\sum_{i=1}^{p-1} \sum_{j=i+1}^{p} |S_{ij}|$.
%consistently with \cite{fan2013large}.
Second, we estimate factor scores and loadings
conditioning on the estimates of $L^{*}$ and $S^{*}$ given by \ref{ours}.

Problem (\ref{ours}) is a least squares one, penalized by a nuclear norm plus $l_1$ norm heuristics,
which has been proved in \cite{fazel2002matrix} to be the tightest convex relaxation of the original NP-hard problem involving $rank(L)$ and $||S||_{0}$.
The optimum is computed via an alternate thresholding algorithm,
composed by a singular value thresholding \citep{cai2010singular}
and a soft-thresholding step \citep{daubechies2004iterative}
(we refer to the supplement of \cite{farne2020large} for more details).
Some variants of (\ref{ours}) have been used to estimate the covariance matrix and its inverse under the low rank plus sparse assumption in \cite{agarwal2012noisy} and \cite{chandrasekaran2012} respectively.

%%ok!
%%simply the eigenvectors of \widehat{L}_{UNALCE}!

Problem (\ref{ours}) can be thought of as an approximate robust PCA problem.
In \cite{farne2020large}, a refined estimation theory for the estimates of $L^{*}$, $S^{*}$ and $\Sigma^{*}$ obtained by (\ref{ours})
is provided assuming the generalized spikiness of the eigenvalues of $L^{*}$ and the generalized element-wise sparsity of $S^{*}$. A characterizing feature of those estimates is that they are both parametrically and algebraically consistent, i.e., the latent rank and the residual sparsity pattern
are exactly recovered.
The effectiveness of problem (\ref{ours}) as a factor model estimation method has been recently studied
in \cite{bai2019rank} as far as parametric consistency is concerned, but no algebraic consistency theory is provided therein. Moreover, the latent eigenvalues must diverge with $p$ in order to ensure parametric consistency.
In this paper we derive finite sample consistency results for factor loadings, factor scores and common components based on the theoretical framework of \cite{farne2020large}, which encompasses a wide range of low rank plus sparse stochastic structures.

The solutions to problem (\ref{ours}) in \cite{farne2020large}
are called $\widehat{\mathit{L}}_{ALCE}$ and $\widehat{\mathit{S}}_{ALCE}$,
where ALCE stands for ALgebraic Covariance Estimator. ALCE estimates are then re-optimized by applying an additional least squares step, leading to the final estimates
$\widehat{\mathit{L}}_{UNALCE}$ and $\widehat{\mathit{S}}_{UNALCE}$ (where UNALCE stands for UNshrunk ALCE).
The main alternative is POET \citep{fan2013large}, a two-step estimator where $L^{*}$ is estimated as the covariance matrix of the top $r$ principal components,
%via PCA as $\widehat{L}_{POET}=\min_{L, rk(L)\leq r} ||\Sigma_n-L||_{F}$
and $S^{*}$ is estimated by soft-thresholding their orthogonal complement.
%of $\Sigma_n$ with respect to $\widehat{L}_{POET}$.
%: $\widehat{S}_{POET}=\min_{S} ||(\Sigma_n-\widehat{L}_{PCA}-S)||_{1}$.
In comparison to \cite{bai2019rank} and \cite{fan2013large},
the estimation framework of this paper gives several advantages:
\begin{enumerate}
\item no need to use any additional criterion to recover the latent rank;
\item intermediately spiked latent eigenvalues are recovered;
\item any residual sparsity pattern is exactly recovered;
\item the sampling theory is relaxed according to the degree of pervasiveness of latent factors
and the degree of sparsity of the residual component;
\item finite sample error bounds are provided.
\end{enumerate}

%In order to start, we will harmonize the assumptions of \cite{fan2013large} and \cite{farne2018finite}
Moving from the assumptions of \cite{bai2003inferential} and \cite{fan2013large}, we now recall the assumptions of \cite{farne2020large} and we introduce new assumptions to establish the consistency of the OLS-based factor scores obtained via (\ref{ours}).

%\bibliographystyle{chicago}%natbib
%\bibliography{factors_bib}
%
%\end{document}

%$\widehat{B}_{ij}=\min_{B_{ij}} ||X_i-B_{ij}f_ij||+$
%% finite fourth moment for u_t and e_t and bounded weak dependence for B_s and C_s.

%%%%%%
%% under consistency! ok!

\section{Factor model estimation under generalized pervasiveness}\label{sec:ols_both}

\subsection{Derivation of estimates}

%%OK! S^{*}=I_p.
%%given the data!

Let us first define the $n\times r$ matrix $F$ as $F'=[f_1 \ldots f_n]$,
the $p\times r$ matrix $B$ as $B'=[b_1 \ldots b_p]$, and
the $n\times p$ data matrix $X$ as $X'=[x_1 \ldots x_n]$.
The factor-model estimates based on the ordinary least squares are derived
as follows: %move from \cite{bai2003inferential}, who
\begin{equation}\label{ols}
\min_{B,F} \frac{1}{pn}\sum_{j=1}^p \sum_{k=1}^n (X_{k,i}-b_j' f_k)^2.
\end{equation}
According to \cite{bai2003inferential}, minimizing (\ref{ols})
amounts to maximizing $tr(F'(XX')F)$. Under the constraints that
$\frac{1}{n} \sum_{k=1}^n \widehat{f}_k \widehat{f}_k'=I_r$ and $\widehat{B}'\widehat{B}$ is diagonal,
(\ref{ols}) is solved by $\widehat{F}_{OLS,1}=\sqrt{n} {U}_n$, where ${U}_n$ is the $n \times r$ matrix of the top $r$ eigenvectors of the $n \times n$ matrix $XX'$, and $\widehat{B}_{OLS,1}'=\frac{1}{n} \widehat{F}_{OLS,1}'X$.

%The asymptotic inferential theory for factor-model estimates obtained by (\ref{ols}) was first derived in \cite{bai2003inferential}.
In \cite{fan2013large}, the asymptotic consistency of the factor-model estimates derived in the same way is proved assuming that the residual covariance matrix is sparse. In particular, uniform asymptotic rates for loadings, factor scores and common components are provided. In this section, we generalize the results of \cite{fan2013large} to a much wider context, assuming the intermediate regimes of latent eigenvalue spikiness and residual element-wise sparsity of \cite{farne2020large}, which encompass the underlying assumptions of \cite{bai2003inferential} and \cite{fan2013large}.

Before proceeding with technicalities, let us explore what happens to factor model estimates imposing alternative constraints to the solutions of (\ref{ols}).
In particular, let us add to $\widehat{B}'\widehat{B}$ diagonal the condition $\sum_{i=1}^p ||\widehat{b}_i||=\max$.
In that case, the solution in $B$ is $\widehat{B}_{OLS,2}={U}_r{\widehat{\Lambda}_r}^{1/2}$, where ${U}_r$ is the $p\times r$ matrix whose columns are the top $r$ eigenvectors of $\Sigma_n$ and $\widehat{\Lambda}_r$ is the diagonal matrix containing the top $r$ eigenvalues of $\Sigma_n$.
Conditionally on $\widehat{B}_{OLS,2}$, the factor scores are then estimated
for $k=1,\ldots,n$ as follows: $\widehat{f}_{k,2}=(\widehat{B}_{OLS,2}'\widehat{B}_{OLS,2})^{-1}\widehat{B}_{OLS,2}'x_k=
\widehat{\Lambda}_r^{-1}\widehat{B}_{OLS,2}'x_k$.

It is worth exploring the relationship between $\widehat{F}_{OLS,1}=\sqrt{n} {U}_n$
and $\widehat{F}_{OLS,2}$, defined as $\widehat{F}_{OLS,2}'=[\widehat{f}_{1,2} \ldots \widehat{f}_{n,2}]$.
Denoting the eigenvalues and the eigenvectors of $X'X/n$ by $\widehat{\lambda}_i$ and ${u}_i$, $i=1,\ldots,p$, we know that the corresponding eigenvalues and eigenvectors of $XX'/n$ are $\widehat{\lambda}_i$
and $\widehat{\lambda}_i^{-1/2}X{u}_i$, respectively.
It follows that $\widehat{F}_{OLS,2}=X U_r\widehat{\Lambda}_r^{-1/2}=\widehat{F}_{OLS,1}/\sqrt{n}$, and $\widehat{F}_{OLS,2}\sqrt{n}=\widehat{F}_{OLS,1}$.
At the same time, we can write $\widehat{B}_{OLS,1}'=1/n \widehat{F}_{OLS,1}'X=1/n \widehat{F}_{OLS,2}'\sqrt{n}X=1/n \widehat{B}_{OLS,2}'X'X \sqrt{n}=
\widehat{B}_{OLS,2}'\Sigma_n \sqrt{n}$.

As a consequence, it follows that any asymptotic rate for $\widehat{B}_{OLS,1}$ and $\widehat{f}_{OLS,1}$
holds for $\widehat{B}_{OLS,2}$ and $\widehat{f}_{OLS,2}$ as well, because the two mapping relationships only depend on $\Sigma_n$, which converges to $\Sigma^{*}$ in relative terms as $p$ diverges under the assumptions of \cite{fan2013large}. This holds for POET-based estimates even under the assumptions of \cite{farne2020large}. %as $p^{\alpha}/\sqrt{n}$ converges.
%It follows that the results for $\widehat{B}_{OLS,2}$ and $\widehat{f}_{OLS,2}$ can be derived by the ones for
%$\widehat{B}_{OLS,1}$ and $\widehat{f}_{OLS,1}$, gathering $\widehat{\lambda}_i$ and $\widehat{\Lambda}$ respectively, which converge to $p^{\alpha}$ as $p^{\alpha}$ diverges.

Considering the estimates $\widehat{B}_{OLS,2}$ and $\widehat{F}_{OLS,2}$ based on ALCE estimates
instead of POET, we note that under the conditions of Corollary 2 in \cite{farne2020large},
i.e. as $p^{\alpha+\delta}/\sqrt{n}$ converges to $0$,
ALCE-based estimates converge to the respective targets.
%%same rate
As a consequence, for a large enough dimension $p$, ALCE and POET estimates are so close to share the relative error bound.

\subsection{Consistency of estimates}

We assume the matrix components $\mf{L}^{*}$ and $\mf{S}^{*}$ to come from the following sets of matrices:
\begin{eqnarray}
\mathcal{B}(r) & = & \{\mf{L} \in \R^{p \times p} \mid {\mf{L}}={\mf{U}\mf{D}\mf{U}^\top}, \mf{U} \in \R^{p
\times r} \mathrm{semi-orthogonal}, \mf{D} \in \R^{r \times r} \mathrm{diagonal}\},\label{var:L}\\
\mathcal{A}(s) & = & \{\mf{S}\in \R^{p\times p} \mid \vert \mathrm{support}(\mf{S})\vert \leq
s\},\label{var:S}
\end{eqnarray}
where $\mathcal{B}(r)$ is the variety of matrices with at most
rank $r$, and
$\mathcal{A}(s)$ is the variety of (element-wise) sparse matrices with %-elementwise
at most $s$ non-zero elements, where $\mathrm{support}(\mf{S})$ is the orthogonal complement of $ker(\mf{S})$
and $\vert \mathrm{support}(\mf{S})\vert$ denotes its dimension.
Denoting by $T(\mf{L}^{*})$ and $\Omega(\mf{S}^{*})$ the tangent spaces to $\mathcal{B}(r)$ and $\mathcal{A}(s)$ respectively, the identifiability of $\mf{L}^{*}$ and $\mf{S}^{*}$ is ensured bounding the following rank-sparsity measures:
\begin{eqnarray}
\xi(T(\mf{L}^{*})) & = & \max_{\mf{M} \in T(\mf{L}^{*}), \Vert\mf{M}\Vert_2 \leq 1} {\Vert\mf{M}\Vert_\infty},\label{xi}\\
\mu(\Omega(\mf{S}^{*})) & = &\max_{\mf{M} \in \Omega(\mf{S}^{*}),\Vert\mf{M}\Vert_\infty \leq 1}\
{\Vert\mf{M}\Vert_2}\label{mu},%\pause \vspace{0.2cm}
\end{eqnarray}
as controlling the product between \ref{xi} and \ref{mu}
ensure that $T(\mf{L}^{*})$ and $\Omega(\mf{S}^{*})$ intersect only at the origin.

%In this Section, we describe consistency results
%for the OLS-based estimates of factor loadings, scores and communalities.
%Our theoretical framework encompasses the settings of \cite{bai2003inferential},
%\cite{fan2013large} and \cite{farne2020large}.
%This happens because consistency results are asymptotic in nature,
%and the three frameworks asymptotically collapse to the same one.
%% a parita' di condizione!
%% qui si fac per factor consistency!
%% prima no perche' ours was proved without!
%% fan VS bai log(p)/n VS \frac{1}{n}!

%to which refer for more detail.
%In this way, we are able to make a founded comparison between the two methods,
%because they are both valid.

%%under the assumptions of \cite{farne2018finite}.
%This happens because any randomness structure entirely depends in both cases
%on the sample covariance matrix under the generalized spikiness and sparsity assumptions.
%In this respect, this section adds value as it describes the performance of
%any method based on principal components and $\Sigma_n$. %POET included.

%We now describe our model assumptions.
We have recalled in the introduction the assumption context of \cite{farne2020large}.
In order to prove our results, we need to recall their six assumptions.
%Assumption 3 bounds the tails of factors and residuals. Assumption 4 controls the degree of sparsity of the sparse component.
%Assumption 5 sets the allowed latent rank and sample size. Assumption 6 ensures the sparsity pattern recovery.
%
%%from that paper.
%%which are reported in the Appendix. %We also borrow the following .
This is needed to allow for intermediate spikiness regimes for latent eigenvalues
and sparsity regimes for the residual component,
to bound the distribution tails of factors and residuals,
 %. Note that the implicit assumption $\delta' \leq 1$ ensures the compatibility with the sparsity assumptions in \cite{fan2013large}.
to impose a prescribed magnitude for the rank and a lower bound for the sample size,
to control for the residual sparsity pattern and to guarantee its recovery.
The assumptions are reported in detail in Appendix \ref{ass:fm20}.

In addition,
the following lower bounds for the smallest latent eigenvalue and the minimum off-diagonal absolute magnitude in the residual component
are crucial for identifiability and recovery of both matrix components.
\begin{Ass}\label{lowerbounds}
\begin{enumerate}
\item The minimum eigenvalue of $\mf{L}^{*}$ ($\lambda_r(\mf{L}^{*})$) is greater than
$C_2 {\psi}/{\xi^2(T)}$.
\item The minimum absolute value of the non-zero off-diagonal entries of $\mf{S}^{*}$,
${S}_{min,off}$, is greater than $C_3\psi/\mu(\Omega)$.
\end{enumerate}
\end{Ass}
Note that Assumption \ref{lowerbounds}.1 ensures both rank recovery and parametric consistency,
while Assumption \ref{lowerbounds}.2 is necessary only to recover the sparse component.
%Assumptions \ref{eigenvalues}-\ref{pr} and \ref{lowerbounds}-1, in fact, ensure the parametric consistency of $\widehat{L}_{ALCE}$ and $\widehat{S}_{ALCE}$ and the exact recovery of $rank(\widehat{L}_{ALCE})$.

We add here a crucial assumption on loadings, residuals and their interaction. %to our theoretical framework.
This assumption generalizes the corresponding assumption of \cite{fan2013large}
to the intermediate spikiness and sparsity regimes.
\begin{Ass}\label{4'}%\cite{POET} Proposition 4'
%[modified]
There exists $M>0$ such that, for all $j\leq p$, $s\leq n$ and $t\leq n$
\begin{enumerate}
\item $||b_j||_{max}<M,$
\item $E[p^{-\alpha/2}(\epsilon_{s}'\epsilon_{t} - E(\epsilon_{s}'\epsilon_{t}))]^4<M$ and
\item $E[||p^{-\alpha/2}\sum_{i=1}^p b_i \epsilon_{t,i}||^4]<M$,
\end{enumerate}
where $\epsilon_{t,i}$ is the $i-$th component of $\epsilon_{t}$.
In addition, $n=o(p^2)$.
\end{Ass}
Assumption \ref{4'} is made weaker \emph{wrt} the corresponding assumption in \cite{fan2013large} according to the true degree of spikiness of latent eigenvalues. Note that we keep the assumption $n=o(p^2)$, in order
to obtain uniform rates for loadings, factor scores, and common components.

%We now add a crucial assumption on loadings, residuals and their interaction. %to our theoretical framework.
%Here, we generalize the version of \cite{fan2013large}, in order to explore what happens under the intermediate spikiness case.
%\begin{Ass}\label{4'}%\cite{POET} Proposition 4'
%%[modified]
%There exists $M>0$ such that, for all $k\leq p$, $i'\leq n$ and $i''\leq n$,
%\begin{enumerate}
%\item $||b_k||_{max}<M,$
%\item $E[p^{-\alpha/2}(\epsilon_i'\epsilon_{i''} - E(\epsilon_i'\epsilon_{i''}))]^4<M$ and
%\item $E[||p^{-\alpha/2}\sum_{i=1}^p b_k \epsilon_{ki'}||^4]<M.$
%\end{enumerate}
%\end{Ass}
%Assumption \ref{4'} is made weaker according to the true degree of spikiness of latent eigenvalues.

%Proof of the main theorem.
%What happens to POET factors in the intermediate spikiness case?

%As previously said, the results we are explaining also hold for POET.

We now focus on factor model estimates.
We follow the inferential framework of \cite{bai2003inferential}, exactly as \cite{fan2013large} does.
We start reasoning on POET factor model estimates based on ordinary least squares. %Recalling that
We define the projection matrix onto the orthogonally rotated true factor space as
$H_{POET}=\frac{1}{n}(\widehat{\Lambda}_r)^{-1}\widehat{F}_{POET}'FB'B$.
%, where $\widehat{\Lambda}_r$ is the diagonal matrix
%composed by the top $r$ sample eigenvalues in descending order.
%$H_{UNALCE}=\frac{1}{n}(\widehat{\Lambda}_{UNALCE})^{-1}\widehat{F}_{UNALCE}'FB'B$.
%\end{slide}
%5-3
%%
%If \textbf{Ass.} \ref{4'} holds in place of \textbf{Ass.} \ref{4_1} and the \textbf{true} $r$ is known,
Then, the following Theorem holds.

%%on the performance of OLS-based factor scores. %Bartlett's factor scores.
%%it is possible to prove that
%\begin{Thm}\label{facPOET}
%Suppose that Assumptions \ref{eigenvalues}, \ref{tails}, \ref{sparsity} and \ref{4'} hold.
%Then, for the OLS factor model estimates based on POET it holds
%$$\max_{j\leq p}||{\widehat{b}_j-H b_j}||=O_p\left(\omega_n\right)$$
%with $\omega_n=p^{\delta/2}\sqrt{\frac{\log p}{n}}+\frac{p^{(1-\alpha)/2}}{p^{\alpha/2}}$ and
%%$p^{1/2-\alpha/4}$
%$$\max_{k\leq n}||{\widehat{f}_k-H f_k}||=O\left({\frac{p^{1-\alpha}}{\sqrt{n}}}+\frac{p^{1-\alpha}n^{1/4}}{p^{\alpha/2}}\right)$$
%and
%$$\max_{j\leq p, k\leq n}||\widehat{b}_j' \widehat{f}_k-b_j'f_k||=O\left(\frac{n^{1/4}p^{1-\alpha}}{p^{\alpha/2}}+\log(n)^{1/b_2}p^{\delta/2}\sqrt{\frac{\log{p}}{n}}\right)$$
%as $p$ and $n$ diverge to infinity.
%%$$\max_{j\leq p, k\leq n}||\widehat{b}_j' \widehat{f}_k-b_j'f_k||=O\left(\frac{p^{1-\alpha} n^{1/4}}{p^{\alpha/2}}\right)+\log(n)^{1/b_2}p^{\alpha/2}\sqrt{\frac{\log{p}}{n}}, \; \mbox{if} \; \alpha\leq 2/3,$$
%%and
%%$$\max_{j\leq p, k\leq n}||\widehat{b}_j' \widehat{f}_k-b_j'f_k||=O\left(\frac{p^{1-\alpha} n^{1/4}}{p^{\alpha/2}}\right)+\log(n)^{1/b_2}p^{1-\alpha}\sqrt{\frac{\log{p}}{n}}, \; \mbox{if} \; 2/3<\alpha \leq 1.$$
%\end{Thm}

%on the performance of OLS-based factor scores. %Bartlett's factor scores.
%it is possible to prove that
\begin{Thm}\label{facPOET}
Suppose that Assumptions \ref{eigenvalues}, \ref{tails}, \ref{sparsity} and \ref{4'} hold.
Then, setting $d=p^{\alpha}$, for the OLS factor model estimates based on POET it holds
$$\max_{j\leq p}\frac{1}{d}||{\widehat{b}_j-H b_j}||=O_p\left(\omega_n\right)$$
with $\omega_n=p^{\alpha+\delta/2}\sqrt{\frac{\log p}{n}}+p$ and
%$p^{1/2-\alpha/4}$
$$\max_{k\leq n}\frac{1}{d}||{\widehat{f}_k-H f_k}||=O\left({\frac{p}{\sqrt{n}}}+\frac{n^{1/4}p}{p^{\alpha/2}}\right)$$
and
$$\max_{j\leq p, k\leq n}\frac{1}{d}||\widehat{b}_j' \widehat{f}_k-b_j'f_k||=O\left(\frac{n^{1/4}p}{p^{\alpha/2}}+\log(n)^{1/b_2}p^{\alpha+\delta/2}\sqrt{\frac{\log{p}}{n}}\right)$$
as $p$ and $n$ diverge to infinity.
%$$\max_{j\leq p, k\leq n}||\widehat{b}_j' \widehat{f}_k-b_j'f_k||=O\left(\frac{p^{1-\alpha} n^{1/4}}{p^{\alpha/2}}\right)+\log(n)^{1/b_2}p^{\alpha/2}\sqrt{\frac{\log{p}}{n}}, \; \mbox{if} \; \alpha\leq 2/3,$$
%and
%$$\max_{j\leq p, k\leq n}||\widehat{b}_j' \widehat{f}_k-b_j'f_k||=O\left(\frac{p^{1-\alpha} n^{1/4}}{p^{\alpha/2}}\right)+\log(n)^{1/b_2}p^{1-\alpha}\sqrt{\frac{\log{p}}{n}}, \; \mbox{if} \; 2/3<\alpha \leq 1.$$
\end{Thm}

Theorem \ref{facPOET} shows that OLS-based POET factor model estimates are still asymptotically consistent under the generalized spikiness and sparsity regimes, provided that the rank $r$ is known.
Otherwise, as reported by Yu and Samworth in the discussion of \cite{fan2013large},
the latent rank may be underestimated
by the information criteria of \cite{bai2002determining} when $\alpha<1$, %used in \cite{fan2013large},
since in that case $\lim_{p,n \rightarrow \infty} P\{IC(r')<P(IC(r))\}>0$, $r'<r$.
%Nevertheless, if $r$ is known or recovered by a consistent method, like ALCE,
%POET theorems still hold, even in the intermediate spikiness case.
Estimated loadings are consistent as long as $\alpha>1/2$ and $n>k_1 p^{\delta}$ for some $k_1>0$.
The consistency of estimated factor scores requires $\alpha > \frac{3}{4}$ and $n=o(p^2)$.
%and $k_1 p^{2(1-\alpha)} < n < k_2 p^{6\alpha-4}$ for some $k_1,k_2>0$.
The consistency of communalities requires both sets of conditions. Note that the asymptotic consistency
requires the convergence condition of $\Sigma_n$ to $\Sigma^{*}$, i.e. the convergence of $p^{\alpha}/\sqrt{n}$ to $0$,
to hold.%
%$\alpha > \frac{3}{4}$ and $n=o(p^2)$.

%%
Concerning ALCE-based factor model estimates, the following result holds.
\begin{Thm}\label{facALCE}
If the assumptions of Theorem \ref{facPOET} and Assumptions \ref{alg}, \ref{pr}, \ref{sc}, and \ref{lowerbounds} hold,
Theorem \ref{facPOET} holds also for the OLS factor model
estimates based on ALCE, setting $d=p^{\alpha+\delta}$.
\end{Thm}
Note that Assumption \ref{pr} encompasses the condition $n>k_1 p^{\delta}$,
and the asymptotic consistency requires the convergence of $\Sigma_n$ to $\Sigma^{*}$,
which holds in this case if $p^{\alpha+\delta}/\sqrt{n} \rightarrow 0$.
Estimated loadings now require the conditions $\alpha+\delta>1/2$
and $n>k_1 p^{\delta}$ for some $k_1>0$ to be consistent,
while the estimated factor scores require $\alpha+\delta > \frac{3}{4}$ and $n=o(p^2)$.
We refer to Appendices \ref{facFan} and \ref{facMe} for the formal proofs.

From the following section, we explore the behaviour of UNALCE-based factor model estimates
whenever the parameters $p$ and $n$ are fixed.
Those estimates, in fact, show very interesting properties as far as numerical stability and fitting properties is concerned.

%Note that the assumptions of \cite{fan2013large} are contained in the ones of \cite{farne2018finite}.
%The additional constraints depend on the need to ensure the exact recovery of rank and sparsity pattern (algebraic consistency).
%$m_p=o(p^\alpha)$.

%% Sketch of the proof

%%displaying!

\section{ALCE and UNALCE in the finite sample}\label{sec:unalce}

In Section \ref{sec:ols_both} we derived the asymptotic consistency of OLS-based factor model estimates obtained via POET and ALCE. In this section, we discuss the optimality properties of factor model estimates based on heuristics (\ref{ours}) when the parameters $p$ and $n$ are fixed. %following heuristics:
%
%\begin{equation}
%\min_{L,S} \frac{1}{2}||(L+S)-\Sigma_n||_{F}^2 +  \lambda
%||L||_{*} + \rho ||S||_{1,off},\label{func:ob1}%
%\end{equation}
%where
%\begin{itemize}
%\item $||S||_{1,off}$ is the $l_1$ norm of
%$S$ excluding the diagonal, i.e. $\sum_{i=1}^p \sum_{j=2}^{p} |s_{ij}|$; %$||S||_0$ has been replaced by
%\item $||L||_{*}$ is the \textbf{nuclear} norm of %$rank (L)$ has been replaced by
%$L=UDU'$, i.e. $\sum_{i=1}^{r} |d_i|=\sum_{i=1}^{r}
%d_i=||diag(D)||_1$.
%\end{itemize}
%We call $\widehat{\mathit{L}}_{ALCE}$ and $\widehat{\mathit{S}}_{ALCE}$ its solutions.
In order to do that, we need to recall two key results of \cite{farne2020large}.
%The first one prescribes that under Assumptions \ref{eigenvalues}-\ref{4'} the rank of L

The first one follows by Theorem \ref{thmMinetop} and Corollary \ref{losses}.
Theorem \ref{thmMinetop} states that the solutions of \ref{ours}
under Assumptions \ref{eigenvalues}-\ref{sc} and Assumption \ref{lowerbounds} are parametrically consistent and
recover exactly the latent rank and the residual sparsity pattern with high probability.
The threshold parameters are set as $\psi=\frac{1}{\xi(T)}\frac{p^{\alpha}}{\sqrt{n}}$,
and $\rho=\gamma \psi$, where $\gamma \in [9\xi(T),1/(6\mu(\Omega))]$.
The resulting estimators are $\widehat{{L}}_{ALCE}$, $\widehat{{S}}_{ALCE}$ and $\widehat{{\Sigma}}_{ALCE}$.
%In particular, Theorem \ref{thmMinetop} states that the latent rank and sparsity pattern
%are exactly recovered with high probability
%for $p$ and $n$ fixed if the prescribed assumptions hold. $\widehat{{L}}$, $\widehat{{S}}$
%and $\widehat{{\Sigma}}=\widehat{{L}}+\widehat{{S}}$ are the ALCE estimates, respectively denoted by
%$\widehat{{L}}_{ALCE}$, $\widehat{{S}}_{ALCE}$ and $\widehat{{\Sigma}}_{ALCE}$.
Corollary \ref{losses} states the finite bounds and the positive definiteness conditions for the residual and the overall ALCE estimates.
Theorem 1 and Corollary \ref{losses} together mean that ALCE estimates are algebraically consistent.

The second key result is related to the finite sample optimization of ALCE estimates.
Let us define ${Y}_{pre}$ and ${Z}_{pre}$
the last updates of the gradient step during the minimization algorithm of (\ref{ours}).
We also define ${\Sigma}_{pre}={Y}_{pre}+{Z}_{pre}$.
In \cite{farne2020large}, it was proved that ALCE estimates can be improved as much as possible
conditioning on ${Y}_{pre}$ and ${Z}_{pre}$. We report here a consequence of that result
relevant for our purposes.

\begin{Thm}\label{mine}
%Suppose that ${{L}}_{pre}$ and ${{S}}_{pre}$ are the last updates in the minimization process of (\ref{ours}).
Suppose that $\widehat{\mathcal{B}}(r)$ and $\widehat{\mathcal{A}}(s)$ are the recovered matrix varieties, and define as $\widehat{{L}}_{ALCE}=\widehat{ U}_{ALCE}\widehat{ D}_{ALCE}\widehat{ U}_{ALCE}'$
the eigenvalue decomposition of $\widehat{{L}}_{ALCE}$.
%\medskip
Assume that ${S}$ has the same off-diagonal elements as $\widehat{{S}}_{ALCE}$
and that the diagonal elements of ${L}+{S}$ are the same as $\widehat{{\Sigma}}_{ALCE}$. %such that its diagonal elements are the same as $\widehat{{\Sigma}}_{ALCE}$ respectively.
Under Assumptions \ref{eigenvalues}-\ref{sc} and Assumption \ref{lowerbounds}, then the minima
\begin{eqnarray}
\min_{{L}\in\widehat{\mathcal{B}}(\widehat{r})}\|{L}-L^{*}\|_{2}\nonumber\\
\min_{{S} \in \widehat{\mathcal{A}}(\widehat{s})}\|{S}-S^{*}\|_{2}\nonumber\\
\min_{{L}\in\widehat{\mathcal{B}}(\widehat{r}), {S} \in \widehat{\mathcal{A}}(\widehat{s})}\|({L}+{S})-\Sigma^{*}\|_{2}\nonumber\\
\min_{{S} \in \widehat{\mathcal{A}}(\widehat{s})}\|{S}^{-1}-S^{*-1}\|_{2}\nonumber\\
\min_{{L}\in\widehat{\mathcal{B}}(\widehat{r}), {S} \in \widehat{\mathcal{A}}(\widehat{s})}\|({L}+{S})^{-1}-\Sigma^{*-1}\|_{2}\nonumber
\end{eqnarray}
conditioning on ${Y}_{pre}$ and ${Z}_{pre}$ are achieved if and only if
$${L}=\widehat{{L}}_{UNALCE}=\widehat{ U}_{ALCE}(\widehat{ D}_{ALCE}+\breve{\psi} {I}_r)\widehat{ U}_{ALCE}'
\quad  \mbox{and}$$
$$\quad
 diag({S})= diag(\widehat{{S}}_{UNALCE})=diag(\widehat{{\Sigma}}_{ALCE})-diag(\widehat{{L}}_{UNALCE})$$
where $\breve{\psi}>0$ is any prescribed threshold parameter. %and UNALCE stands for UNshrunk ALCE. %in (\ref{ours}).
\end{Thm}
Theorem \ref{mine} %holds because if we re-add the threshold $\psi$ to the eigenvalues of $\widehat{L}_{ALCE}$ we obtain $\widehat{L}_{UNALCE}$, and that process minimizes the loss from the target $L^{*}$ conditioning on $Y_{pre}$ into the recovered low rank variety.
%The same holds for $\widehat{S}_{UNALCE}$ into the recovered sparse variety and for the overall covariance estimator $\widehat{\Sigma}_{UNALCE}=\widehat{L}_{UNALCE}+\widehat{S}_{UNALCE}$, keeping fixed the diagonal elements of $\widehat{\Sigma}_{ALCE}$ and the off-diagonal sparsity pattern of $\widehat{S}_{UNALCE}$.
%Therefore, conditioning on ${Y}_{pre}$ and ${Z}_{pre}$,
states that the UNALCE estimates of $L^{*}$,$S^{*}$,$\Sigma^{*}$,$S^{*-1}$,$\Sigma^{*-1}$ show the least possible errors in spectral (and Frobenius) norm
within the class of algebraically consistent estimates,
conditioning on the data.
%Consequently, they show systematic gains over ALCE estimates.
%$\widehat{\Sigma}_{UNALCE}$ also shows the least possible loss from the sample covariance matrix $\Sigma_n$.
%The conditioning properties of $\widehat{L}_{UNALCE}$, $\widehat{S}_{UNALCE}$, $\widehat{\Sigma}_{UNALCE}$ are also reshaped
%with respect to their ALCE counterparts. %we refer to \cite{farne2018finite}.
We note that Weyl's theorem ensures that the absolute errors of
UNALCE individual eigenvalues also have the minimum possible upper bound
under the same assumptions.
We refer to Appendix \ref{proofs_3} for the proofs.
%$|\widehat{\lambda}_{L,i}-\lambda_{L,i}|\leq ||\widehat{L}-L^*||$, $i=1,\ldots,r$,
%and $|\widehat{\lambda}_{\Sigma,i}-\lambda_{\Sigma,i}|\leq ||\widehat{\Sigma}-\Sigma^*||$,
%which means
%We refer to \cite{farne2020large} for the details.

%%

\section{Optimality properties of UNALCE estimates}\label{sec:eig_conc}

%%Further optimality properties of UNALCE estimates (given the finite sample)

We now analyze the parametric and algebraic properties of $(\widehat{L}_{UNALCE},\widehat{S}_{UNALCE})$
with respect to $(\widehat{L}_{ALCE},\widehat{S}_{ALCE})$ and $(\widehat{L}_{POET},\widehat{S}_{POET})$,
and their impact on factor model estimates.
Proving the consistency of the estimates obtained by (\ref{ours}) involves sub-differential methods and fixed point theorems.
The reference norm to assess consistency is the dual norm of the cartesian space $\mathcal{Y}=\mathcal{B}(r)\oplus\mathcal{A}(s)$, which is
$g_{\gamma}(\widehat{S}-S^{*},\widehat{L}-L^{*})=\max\left({||\widehat{L}-L^{*}||_2,\frac{||\widehat{S}-S^{*}||_{\infty}}{\gamma}}\right).$
In \cite{luo2011high}, it is shown that the proof requires to solve three algebraic problems.
The first one requires the minimization of (\ref{ours}) under the constraint $(L,S)\in \mathcal{M}$,
where $\mathcal{M}$ is the class of low rank matrices $L$ and sparse matrices $S$ satisfying the following conditions
$$||\mathbb{P}_{T'^\perp} (L-L^*)||\leq \xi(T)\psi,$$
$$g_{\gamma}(\Sigma-\Sigma^{*},\Sigma-\Sigma^{*})\leq 11 \psi,$$
%($\psi$ is the probabilistic bound $\frac{1}{\xi(T)}\frac{p}{\sqrt{n}}$ and
provided that $\Sigma=L+S$, $\mathbb{P}$ is the projection operator and $T'$ is a manifold sufficiently close to the tangent space $T$.
As a consequence, those constraints hold for $\widehat{L}_{ALCE}$, $\widehat{S}_{ALCE}$, and $\widehat{\Sigma}_{ALCE}=\widehat{L}_{ALCE}+\widehat{S}_{ALCE}$.
From this consideration, we can derive the following corollary.

%Assuming that the off-diagonal elements of $S$ are the same as $\widehat{{S}}_{ALCE}$
%and the diagonal elements of $\Sigma$ are the same as $\widehat{{\Sigma}}_{ALCE}$, we obtain the solutions
%$L=\widehat{L}_{UNALCE}$, $S=\widehat{S}_{UNALCE}$. %Note that the minimum in $L$ is $0$, while the one in $S$ is not.
%Therefore, if we write $\min_{\Sigma \in \mathcal{Y}} ||{\Sigma}-\Sigma^*||_{F}^2=||{\Sigma}-\widehat{\Sigma}_{ALCE}||_{F}^2+||\widehat{\Sigma}_{ALCE}-\Sigma^*||_{F}^2$,
%i.e. the MSE of $\widehat{\Sigma}$ is minimum conditioning on $\widehat{\Sigma}_{ALCE}$ because $||{\Sigma}-\widehat{\Sigma}_{ALCE}||_{F}^2$ is minimum.
%The same can be written for $L$ and $S$, where $||{L}-\widehat{L}_{ALCE}||_{F}^2=0$ and $||{S}-\widehat{S}_{ALCE}||_{F}^2 \ne 0$.

\begin{Coroll}\label{Coroll_pl}
In general, it holds
$$||\mathbb{P}_{T'^\perp} (\widehat{L}_{UNALCE}-L^*)||\leq (C+1)\psi$$
$$g_\gamma(\widehat{\Sigma}_{UNALCE}-\Sigma^{*},\widehat{\Sigma}_{UNALCE}-\Sigma^{*})\leq (C+2) \psi,$$
where $C$ is the positive constant of Theorem \ref{thmMinetop}.

Conditionally on $Y_{pre}$ and $Z_{pre}$, it holds
$$||\mathbb{P}_{T'^\perp} (\widehat{L}_{ALCE}-L^*)||-||\mathbb{P}_{T'^\perp} (\widehat{L}_{UNALCE}-L^*)||
\leq \psi,$$
$$g_\gamma(\widehat{\Sigma}_{ALCE}-\Sigma^{*},\widehat{\Sigma}_{ALCE}-\Sigma^{*})-
g_\gamma(\widehat{\Sigma}_{UNALCE}-\Sigma^{*},\widehat{\Sigma}_{UNALCE}-\Sigma^{*})\leq \psi.$$

%$0 \leq c_1 \leq \psi$ and  $0 \leq c_2 \leq \psi$.
\end{Coroll}
We refer to Appendix \ref{diag_trace} for a discussion of the algebraic and parametric properties
of POET and UNALCE component error estimates.
%and the optimality properties of the traces and diagonals of the
%The condition $g_\gamma(\mathcal{P}_{\mathcal{T}^\perp}(\widehat{S}+\widehat{L}-\Sigma_n)) < \psi$,
%which is sufficient to ensure that (\ref{ours}) is minimum (see \cite{luo2011high}),
%still holds for UNALCE, because it still holds $\mathcal{P}_{\mathcal{Y}^\perp} A^\perp[\widehat{L}+\widehat{S}-\Sigma_n]=Z=(-\psi\gamma sign(S^*),-\psi UU')$,
%since the projections onto the orthogonal complements of ${T}$ and ${\Omega}$ are unvaried.
%and $$ $L+S-\Sigma_n$ is unaltered,
%due to the fact that $\widehat{L}+\widehat{S}-\Sigma_n$ is equal for UNALCE and ALCE
%(re-optimized least squares!)

%\section{UNALCE optimal eigenvalue concentration}

We now report a crucial property of the eigenvalues of UNALCE estimates.
\begin{Thm}\label{eigen}
%the eigenvalues of $L^{*}$,$S^{*}$,$\Sigma^{*}$,$S^{*-1}$, $\Sigma^{*-1}$ %estimated by UNALCE show the minimum variance %possible, i.e. they are the most concentrated possible around their true means
Let us define $\mu_{L}=tr(L^{*})/p$, $\mu_{S}=tr(S^{*})/p$, $\mu_{\Sigma}=tr(\Sigma^{*})/p$, $\mu_{S^{-1}}=tr(S^{*-1})/p$, $\mu_{\Sigma^{*-1}}=tr(\Sigma^{*-1})/p$.
Under the assumptions of Theorem \ref{mine}, the following statements hold:
%$\mu_{S}$, $\mu_{\Sigma}$, $\mu_{S^{-1}}$,$ \mu_{\Sigma^{-1}}$ given the sample covariance matrix $\Sigma_n$:
\begin{eqnarray}
\widehat{L}_{UNALCE}&=&\min_{L \in \widehat{\mathcal{B}}(\widehat{r})}\frac{1}{p} E\left[\sum_{i=1}^p (\widehat{\lambda}_{L,i}-\mu_{L})^2\vert \Sigma_n\right],\nonumber\\
\widehat{S}_{UNALCE}&=&\min_{S \in \widehat{\mathcal{A}}(\widehat{s})}\frac{1}{p} E\left[\sum_{i=1}^p (\widehat{\lambda}_{S,i}-\mu_{S})^2\vert \Sigma_n\right],\nonumber\\
\widehat{\Sigma}_{UNALCE}&=&\min_{\Sigma \in \widehat{\mathcal{Y}}}\frac{1}{p} E\left[\sum_{i=1}^p (\widehat{\lambda}_{\Sigma,i}-\mu_{\Sigma})^2\vert \Sigma_n\right],\nonumber\\
\widehat{S}_{UNALCE}^{-1}&=&\min_{S \in \widehat{\mathcal{A}}(\widehat{s})}\frac{1}{p} E\left[\sum_{i=1}^p (\widehat{\lambda}_{S^{-1},i}-\mu_{S^{-1}})^2 \vert \Sigma_n \right], \nonumber\\
\widehat{\Sigma}_{UNALCE}^{-1}&=&\min_{\Sigma \in \widehat{Y}} \frac{1}{p} E\left[\sum_{i=1}^p(\widehat{\lambda}_{\Sigma^{-1},i}-\mu_{\Sigma^{-1}})^2\vert \Sigma_n \right] \nonumber.
\end{eqnarray}
\end{Thm}

Theorem \ref{eigen} states that the expected dispersion of UNALCE estimated eigenvalues around
the true mean eigenvalue is the minimum possible within the classes of algebraically consistent estimates,
thus outperforming both ALCE and POET. This important result follows from the eigenvalue dispersion lemma
of \cite{ledoit2004well} (see Appendix \ref{proofs_4} for the proof).

%another relevant result on the eigen-structure of our estimates.
According to \cite{bun2017cleaning}, we can define the empirical spectral density function (ESD)
of a matrix $M$, $\rho(z)_M^p$, $z \in R^{+}$, as follows: $\rho(z)_{M}^p=\frac{1}{p}\sum_{i=1}^p \delta(z-\lambda_{M,i})$,
where $\delta(z-\lambda_{M,i})$ is the Dirac-delta function.
We know that the $k-$th moment of the ESD of $M$ is equal to $p^{-1}tr(M^k)$.
The limit of $\rho(z)_M^p$ as $p$ and $z$ go to infinity, that is the limiting spectral distribution (LSD),
is defined as $\rho(z)_M=\lim_{p\rightarrow \infty} \rho(z)_M^p$.
%The Stieltjes transform $g(z)_M$ of $\rho(z)_M$ is defined as
%$g(z)_M=\int\frac{\rho(u)_M}{z-u}du$, and

Given these definitions, from Theorem \ref{eigen} we can state Corollary \ref{esd}.
\begin{Coroll}\label{esd}
Under the assumptions of Theorems \ref{mine},
the second moments of $\rho(z)_{\widehat{L}_{UNALCE}-L^{*}}^p$, $\rho(z)_{\widehat{S}_{UNALCE}-S^{*}}^p$, $\rho(z)_{\widehat{\Sigma}_{UNALCE}-\Sigma^{*}}^p$, $\rho(z)_{\widehat{S}^{-1}_{UNALCE}-S^{*-1}}^p$, $\rho(z)_{\widehat{\Sigma}^{-1}_{UNALCE}-\Sigma^{*-1}}^p$ are the minimum possible within the classes of algebraically consistent estimates.
As $\frac{p^{\alpha+\delta}}{\sqrt{n}} \rightarrow 0$, %and $p\rightarrow \infty$,
the first moments of $\rho(z)_{\widehat{L}_{UNALCE}-L^{*}}$, $\rho(z)_{\widehat{S}_{UNALCE}-S^{*}}$, $\rho(z)_{\widehat{\Sigma}_{UNALCE}-\Sigma^{*}}$, $\rho(z)_{\widehat{S}^{-1}_{UNALCE}-S^{*-1}}$, $\rho(z)_{\widehat{\Sigma}^{-1}_{UNALCE}-\Sigma^{*-1}}$ converge to zero.
%and the second moments of the same Stieltjes transforms
\end{Coroll}
Corollary \ref{esd}, proved in Appendix \ref{proof_esd},
states that target eigenvalues are estimated in the best possible way by UNALCE
within the classes of algebraically consistent estimates.

\section{Bartlett's and Thompson's factor scores optimality}\label{sec:bt_fac}

%%conditional variance of the eigenvalue product around their grand mean.
In this section, we prove that Bartlett's and Thompson's factor scores estimates based on UNALCE
show the minimum loss given the finite sample.
First, we state the optimality of the UNALCE-based loading matrix, $\widehat{B}_{UNALCE}=U_{ALCE}\sqrt{D_{UNALCE}}$, with $\widehat{D}_{UNALCE}=\widehat{D}_{ALCE}+\breve{\psi}I_r$, by the following Corollary.
\begin{Coroll}\label{Coroll_B}
Under the assumptions of Theorem \ref{mine} and Assumption \ref{4'}, the constraints $\widehat{B}'\widehat{B}$ diagonal and $\sum_{i=1}^p ||\widehat{b}_i||=\max$, the minimum $$\min_{\widehat{B},\widehat{L}=\widehat{B}\widehat{B}'\in\widehat{\mathcal{B}}(\widehat{r})}||\widehat{B}-B||$$
is for $\widehat{B}=\widehat{{B}}_{UNALCE}$.
\end{Coroll}
Corollary \ref{Coroll_B} is a direct consequence of Theorem \ref{eigen}.

Then, we define Bartlett's factor scores estimates for the observation $k$, $k=1,\ldots,n$, as follows: $\widehat{f}_{k,B}=(\widehat{B}'\widehat{S}^{-1}\widehat{B})^{-1}\widehat{B}'\widehat{S}^{-1}x_k$.
They simply are the GLS estimates of factor scores conditioning on the data.
We can also define the projections onto the estimated latent space, also called communalities,
as $\widehat{B}\widehat{f}_{k,B}$ for $k=1,\ldots,n$.
The true Bartlett's factors are defined as $f_{k,B}=({B}'{S}^{*-1}{B})^{-1}{B}'{S}^{*-1}x_k$.
The following result for Bartlett's factor scores and projections
onto the latent space based on UNALCE holds.
\begin{Thm}\label{bartlett_opt}
Under the assumptions of Theorem \ref{mine} and Assumption \ref{4'}, the minima for $k=1,\ldots,n$
\begin{eqnarray}
\min_{\widehat{B},\widehat{L}=\widehat{B}\widehat{B}'\in\widehat{\mathcal{B}}(\widehat{r}),\widehat{S} \in \widehat{\mathcal{A}}(\widehat{s})}||\widehat{f}_{k,B}-f_{k,B}||\nonumber\\
\min_{\widehat{B},\widehat{L}=\widehat{B}\widehat{B}'\in\widehat{\mathcal{B}}(\widehat{r}),\widehat{S}
\in \widehat{\mathcal{A}}(\widehat{s})}
||\widehat{B}\widehat{f}_{k,B}-B f_{k,B}||
\end{eqnarray}
conditioning on ${Y}_{pre}$ and ${Z}_{pre}$ are achieved if and only if
$\widehat{B}=\widehat{{B}}_{UNALCE}$ and
$\widehat{S}=\widehat{{S}}_{UNALCE}$.
\end{Thm}
Theorem \ref{bartlett_opt} states that Bartlett's factor scores and communalities estimated by UNALCE
are the most precise given the finite sample within the classes of algebraically consistent estimates for $B$ and $S^{*}$.

Suppose now that the bivariate distribution $(x_k,f_k)$,
$k=1,\ldots,n$, is normal, i.e.
$$\left(\begin{array}{rl}x_k\\ f_k \end{array} \right)\sim NMV\left[\left(\begin{array}{rl}\mu\\
0 \end{array}\right),\left(\begin{array}{rl}BB'+S^{*} & L^{*}\\ L^{*} & I_r\end{array}\right)\right].$$
As a consequence, from the Bayesian point of view, we can derive the following a posteriori expected value for $f_k$:
$$E(f_k|x_k)=B'(BB'+S^{*})^{-1}x_k.$$
%where $\bar{x}$ is the sample mean.
%\pause\vspace{0.2cm}
Thompson's estimates of factor scores
are the estimates of this expected value:
$\widehat{f}_{k,T}=\widehat{B}'(\widehat{B}\widehat{B}'+\widehat{S})^{-1}x_k$.
The corresponding Thompson's true factor scores are $f_{k,T}=B'(BB'+S)^{-1}x_k$.
The following Theorem on the performance of Thompson's estimates of factor scores and communalities
based on UNALCE holds.

\begin{Thm}\label{thompson_opt}
Under the assumptions of Theorem \ref{mine} and Assumption \ref{4'}, the minima for $k=1,\ldots,n$
\begin{eqnarray}
\min_{B,{L}={B}{B}'\in\widehat{\mathcal{B}}(\widehat{r}),{S} \in \widehat{\mathcal{A}}(\widehat{s})}||\widehat{f}_{k,T}-f_{k,T}||\nonumber\\
\min_{B,{L}={B}{B}'\in\widehat{\mathcal{B}}(\widehat{r}),{S} \in \widehat{\mathcal{A}}(\widehat{s})}
||\widehat{B}\widehat{f}_{k,T}-B f_{k,T}||
\end{eqnarray}
are achieved if and only if
$\widehat{B}=\widehat{{B}}_{UNALCE}$ and
$\widehat{S}=\widehat{{S}}_{UNALCE}$.
\end{Thm}
Theorem \ref{thompson_opt} states the same optimality properties of Theorem \ref{bartlett_opt} for Thompson's estimates.
Both proofs, reported in Appendices \ref{proof_bart} and \ref{proof_thom}, rely on Theorem \ref{eigen} and Corollary \ref{esd},
and involve results on the inverse of a matrix sum.
We stress that the optimality of UNALCE with respect to the estimates of factor scores holds both against ALCE and POET, as long as the asymptotic rates of Theorem \ref{facPOET} converge to $0$.

Bartlett's and Thompson's estimates based on UNALCE converge to $\widehat{f}_{OLS,UNALCE}$, because as $p$ and $n$ diverge respecting the condition $p^{\alpha+\delta}/\sqrt{n} \rightarrow 0$, $\widehat{S}$ converges to $I_p$ in the former case, and $I_p$ is negligible with respect to $\widehat{B}\widehat{B}'$ in the second case.
Therefore, the uniform rates derived in Section \ref{sec:ols_both} asymptotically hold for UNALCE Bartlett's and Thompson's estimates too.

\section{Simulation study}\label{sec:sim}

\subsection{Simulation settings}

In this section, we test the validity of Theorems \ref{bartlett_opt} and \ref{thompson_opt}
on some data simulated for that purpose. Here we report our main simulation parameters:
\begin{enumerate}
\item the dimension $p$, the sample size $n$;
\item the rank $r$ and the condition number $cond({L}^*)=\lambda_{max}({L}^*)/\lambda_{min}({L}^*)$
of the low rank component ${L}^*$;% in spectral norm;
\item the trace of ${L}^*$, $\tau \theta p$, where $\tau$ is a magnitude parameter
and $\theta$ is the proportion of variance explained by ${L}^*$;
\item the number of off-diagonal non-zeros $s$ in the sparse component ${S}^*$;
\item the proportion of non-zeros $\pi_s$ over the number of off-diagonal elements;
\item the proportion of the (absolute) residual covariance $\rho_{{S}^*}$;
\item $N=100$ replicates for each setting.
\end{enumerate}
Essentially, the low rank component is simulated by setting $r$ equispaced eigenvalues with sum $\tau \theta p$ and deriving an orthonormal $r-$dimensional basis by Gram-Schmidt algorithm.
The residual variances are simulated by a $p-$dimensional Dirichlet distribution with sum $1-\tau \theta p$, and then matched to the previously simulated diagonal elements of the low rank component according to their relative magnitude.
The off-diagonal elements are first simulated entry-wise by exploiting Cauchy-Schwartz inequality.
The smallest $p(p-1)/2-s$ absolute off-diagonal elements are then set to $0$.
The detailed simulation algorithm is reported in \cite{farne2016large}.

The main parameters of simulated settings are reported in Tables \ref{sett} and \ref{specond}.
We can see that Setting 1 presents not so spiked eigenvalues and a very sparse residual component. This is the most consistent setting with UNALCE assumptions. Setting 2 has spiked eigenvalues and a far less sparse residual. Settings 3 and 4 are intermediately spiked and sparse but present a much lower $p/n$ ratio. In particular, while Settings 1 and 2 have $p/n=10$, Setting 3 has $p/n=1$ and Setting 4 has $p/n=0.5$. Setting 4 is the most consistent with POET assumptions.

In each setting, the eigenvalues of ${L}^*$ and ${\Sigma}^*$ almost overlap, while
the eigenvalues of ${S}^*$ are much smaller. Note that the minimum allowed off-diagonal residual element in absolute value, ${S}_{min,off}$, decreases from Setting 1 to Setting 4.
%\subsection{Simulated settings}

\begin{table}
  \caption{\label{sett} Simulated settings: parameters}
  \centering
  %%\hline
  \fbox{
  \begin{tabular}{cccccccccccc}
    %\hline
    % after \\: %\hline or \cline{col1-col2} \cline{col3-col4} ...
    \textbf{Setting}  &  $p$ &  $n$ &  $p/n$ &  $r$ &  $\theta$ &  $c$ &  $\pi_s$ &  $\rho_{\mathbf{S}^*}$&  \mbox{\textbf{spikiness}}&  \mbox{\textbf{sparsity}}\\% &  $\VertL\Vert$\\
    %\hline
    \textbf{1} & $100$ & $1000$ & $0.1$ & $4$ & $0.7$ & $2$ & $0.0238$ & $0.0045$ & \mbox{low} &\mbox{high}\\% & $23.33$\\
    %%\hline
    %  & $p$ & $n$ & $r$ & $\tau$ & $\theta$ & $c$ & $s$ & $\rho_{\widehat{S}}$\\
    %%\hline
    %%\hline
    % after \\: %\hline or \cline{col1-col2} \cline{col3-col4} ...
    %& $p$ & $n$ & $r$ & $\tau$ & $\theta$ & $c$ & $s$ & $\rho_{\widehat{S}}$\\
    %%\hline
    %\textbf{2} & $100$ & $1000$ & $0.1$ & $4$ & $0.7$ & $4$ & $0.0677$ & $0.0048$ & \mbox{middle} & \mbox{middle}\\% & $28$\\
    \textbf{2} & $100$ & $1000$ & $0.1$ & $3$ & $0.8$ & $4$ & $0.1172$ & $0.0072$ & \mbox{high} & \mbox{low}\\% & $128$\\
    %%\hline
    \textbf{3} & $150$ & $150$ & $1$ & $5$ & $0.8$ & $2$ & $0.0320$ & $0.0033$ & \mbox{middle} & \mbox{middle}\\% & $32$\\
    \textbf{4} & $200$ & $100$ & $2$ & $6$ & $0.8$ & $2$ & $0.0366$ & $0.0039$ & \mbox{middle} & \mbox{middle}\\% & $35.56$\\
    %\hline%yes!
  \end{tabular}
}
\end{table}

\begin{table}
  \caption{\label{specond} Simulated settings: spectral norms and condition numbers}
  \centering
  \fbox{
  \begin{tabular}{ccccccccc}
    % after \\: %\hline or \cline{col1-col2} \cline{col3-col4} ...
    % &  $p$ &  $n$ &  $r$ &  $\tau$ &  $\theta$ &  $c$ &  $s$ &  $\rho_{\widehat{S}}$\\
    %\hline
    \textbf{Setting} &  $\Vert{L}^*\Vert_{2}$ &  $\lambda_r({L}^*)$ &  $cond({{L}^*})$ &  $\Vert{S}^*\Vert_{2}$ &  ${S}_{min,off}$ &  $cond({{S}^*})$ &  $\Vert{\Sigma}^*\Vert_{2}$ &  $cond({{\Sigma}^*})$\\
    %\hline
    \textbf{1} & $23.33$ & $11.67$ & $2$ & $3.78$ & $0.0275$ & $2.26e+07$ & $24.49$ & $9.49e+07$\\
    %%\hline
    %  & $p$ & $n$ & $r$ & $\tau$ & $\theta$ & $c$ & $s$ & $\rho_{\widehat{S}}$\\
    %%\hline
    %%\hline
    % after \\: %\hline or \cline{col1-col2} \cline{col3-col4} ...
    %& $p$ & $n$ & $r$ & $\tau$ & $\theta$ & $c$ & $s$ & $\rho_{\widehat{S}}$\\
    %%\hline
    %\textbf{2} & 28 & 7 & 4 & 2.57 & 0.0143 & 3.80e+07 & 28.83 & 4.04e+07\\
    \textbf{2} & $128$ & $32$ & $4$ & $5.58$ & $0.0226$ & $2.53e+05$ & $130.14$ & $4.07e+06$\\
    %%\hline
    \textbf{3} & $32$ & $16$ & $2$ & $2.56$ & $0.0161$ & $2.35e+13$ & $32.48$ & $1.58e+10$\\
    \textbf{4} & $35.56$ & $17.78$ & $2$ & $4.69$ & $0.0138$ & $1.17e+13$ & $36.39$ & $3.09e+09$\\
    %\hline
  \end{tabular}
}
\end{table}

%Suppose that Assumptions \ref{eigenvalues}, \ref{tails}, \ref{sparsity} and \ref{4'} hold.
%Then, for the OLS factor model estimates based on POET it holds
For each setting, and each of the $h=1,\ldots,100$ replicates,
we simulate $n$ data vectors $z_{h,k}$, $k=1,\ldots,n$, and
we define the respective unbiased sample covariance matrix as $\Sigma_{n,h}$,
the respective spectral decomposition as $\widehat{U}_{h}\widehat{\Lambda}_{h}\widehat{U}_{h}'$,
the sample covariance matrix based on the top $r$ principal components as $\widehat{U}_{h}\widehat{\Lambda}_{h,r}\widehat{U}_{h}'$.
We then apply the minimization algorithm of (\ref{ours}) on $\Sigma_{n,h}$ to get ALCE estimates,
and we derive the subsequent UNALCE and POET covariance estimates: %as follows:
$\widehat{L}_{h,UN}=\widehat{B}_{h,UN}\widehat{B}_{h,UN}'=
(\widehat{U}_{h,UN}\widehat{\Lambda}_{h,UN}^{1/2})(\widehat{U}_{h,UN}\widehat{\Lambda}_{h,UN}^{1/2})'$,
%\vspace{0.2cm}
$\widehat{L}_{h,P}=\widehat{\Sigma}_{h,r}=\widehat{B}_{h,P}\widehat{B}_{h,P}'=
(\widehat{U}_{h,r}\widehat{\Lambda}_{h,r}^{1/2})(\widehat{U}_{h,r}\widehat{\Lambda}_{h,r}^{1/2})'$,
$\widehat{S}_{h,UN}$, $\widehat{S}_{h,POET}$.

Consequently, we derive
Bartlett's estimates ($k=1,\ldots,n$):
%\vspace{0.2cm}
$$\widehat{f}_{h,i,UNALCE,Bartlett}
%\vspace{0.1cm}
=(\widehat{B}_{h,UN}'(\widehat{S}_{h,UN})^{-1}\widehat{B}_{h,UN})^{-1}\widehat{B}_{h,UN}'(\widehat{S}_{h,UN})^{-1}(z_{h,k}-\bar{z}_h),$$
%\vspace{0.2cm}
$$\widehat{f}_{h,i,POET,Bartlett}%=$
%\vspace{0.2cm}
=(\widehat{B}_{h,P}'(\widehat{S}_{h,P})^{-1}\widehat{B}_{h,P})^{-1}\widehat{B}_{h,P}'(\widehat{S}_{h,P})^{-1}(z_{h,k}-\bar{z}_h),$$
and Thompson's estimates of factor scores:
$$\widehat{f}_{h,i,UNALCE,Thompson}
%\vspace{0.1cm}
=\widehat{B}_{h,UN}'(\widehat{\Sigma}_{h,UN})^{-1}(z_{h,k}-\bar{z}_h),$$
%\vspace{0.2cm}
$$\widehat{f}_{i,POET,Thompson}%=$
%\vspace{0.2cm}
=\widehat{B}_{h,P}'(\widehat{\Sigma}_{h,P})^{-1}(z_{h,k}-\bar{z}_h).$$

Defining $H=\frac{1}{n}\widehat{\Lambda}_r^{-1}\widehat{F}'FB'B$,
we calculate the metrics of Theorem \ref{facPOET} for both POET and UNALCE
Bartlett's and Thompson's estimates and
for each replicate $h=1,\ldots,100$ :
$$Loss_B(h)=\max_{j\leq p}||{\widehat{b}_{h,j}-H b_j}||,$$
%$p^{1/2-\alpha/4}$
$$Loss_f(h)=\max_{k\leq n}||{\widehat{f}_{h,k}-H f_k}||,$$
and
$$Loss_{Bf}(h)=\max_{j\leq p, k\leq n}||\widehat{b}_{h,j}' \widehat{f}_{k,j}-b_j'f_k||.$$
%$$\max_{j\leq p, k\leq n}||\widehat{b}_j' \widehat{f}_k-b_j'f_k||=O\left(\frac{p^{1-\alpha} n^{1/4}}{p^{\alpha/2}}\right)+\log(n)^{1/b_2}p^{\alpha/2}\sqrt{\frac{\log{p}}{n}}, \; \mbox{if} \; \alpha\leq 2/3,$$
%and
%$$\max_{j\leq p, k\leq n}||\widehat{b}_j' \widehat{f}_k-b_j'f_k||=O\left(\frac{p^{1-\alpha} n^{1/4}}{p^{\alpha/2}}\right)+\log(n)^{1/b_2}p^{1-\alpha}\sqrt{\frac{\log{p}}{n}}, \; \mbox{if} \; 2/3<\alpha \leq 1.$$
In addition, we calculate the projection of the low rank error matrix onto the orthogonal complement of $L^{*}$
and we measure for each replicate $h$ the magnitude of that matrix in spectral norm for POET and UNALCE:
$$PrErr_{h,P}=||\mathbb{P}_L (\widehat{L}_{h,P}-L^{*})||,$$
$$PrErr_{h,UN}=||\mathbb{P}_L (\widehat{L}_{h,UN}-L^{*})||.$$

Finally, we calculate the means, variances, medians and median absolute deviations
of $Loss_B$, $Loss_f$, $Loss_{Bf}$ and $PrErr$ over the $N$
replicates, both for UNALCE and POET.

\subsection{Simulation results}

\begin{table}[h]
        \centering
        \caption{Simulation results: means and standard deviations of the four sample losses calculated for Bartlett's factor scores over 100 runs.}
        \label{res_all}
        %\resizebox{8cm}{3cm} {
        \begin{tabular}{|c|c|c|c|c|c|c|c|c|c|c|c|c|c|c|}
            \cline{1-10}
             & & \multicolumn{2}{|c|}{Setting 1} & \multicolumn{2}{|c|}{Setting 2} & \multicolumn{2}{|c|}{Setting 3} & \multicolumn{2}{|c|}{Setting 4}\\
            \cline{3-10}
            & & UNALCE & POET & UNALCE & POET & UNALCE & POET & UNALCE & POET\\
            \hline
            $Loss_B$ & mean & 2.8385 & 3.186 & 4.5077 & 4.701 & 3.5768 & 3.773 & 4.5555 & 4.8756\\
                & std & 0.1045 & 0.1586 & 0.1407 & 0.1829 & 0.2104 & 0.2564 & 0.4084 & 0.5361\\
            $Loss_f$ & mean & 0.1928 & 0.3566 & 0.2478 & 0.2916 & 0.3371 & 0.3848 & 0.4926 & 0.5305\\
                & std & 0.0266 & 0.0438 & 0.0796 & 0.0344 & 0.0632 & 0.073 & 0.1105 & 0.1003\\
            $Loss_{Bf}$ & mean & 0.9652 & 2.0299 & 2.1791 & 2.6577 & 2.1976 & 2.424 & 3.1832 & 3.5572\\
                & std & 0.1177 & 0.2435 & 0.5464 & 0.2565 & 0.2227 & 0.2749 & 0.3805 & 0.4871\\
            $PrErr$ & mean & 0.9064 & 1.921 & 2.674 & 3.2001 & 2.8525 & 3.1922 & 4.4129 & 5.0542\\
                & std & 0.1192 & 0.2277 & 0.2927 & 0.4532 & 0.3206 & 0.395 & 0.4614 & 0.7066\\
            \hline
        \end{tabular}
        %}
\end{table}

\begin{table}[h]
        \centering
        \caption{Simulation results: medians and median absolute deviations of the four sample losses calculated for Bartlett's factor scores over 100 runs.}
        \label{res_all2}
        %\resizebox{8cm}{3cm} {
        \begin{tabular}{|c|c|c|c|c|c|c|c|c|c|c|c|c|c|c|}
            \cline{1-10}
             & & \multicolumn{2}{|c|}{Setting 1} & \multicolumn{2}{|c|}{Setting 2} & \multicolumn{2}{|c|}{Setting 3} & \multicolumn{2}{|c|}{Setting 4}\\
            \cline{3-10}
            & & UNALCE & POET & UNALCE & POET & UNALCE & POET & UNALCE & POET\\
            \hline
            $Loss_B$ & median & 2.848 & 3.1894 & 4.4935 & 4.6756 & 3.5614 & 3.7427 & 4.469 & 4.7674\\
                & mad & 0.0851 & 0.1258 & 0.1125 & 0.1492 & 0.1703 & 0.2085 & 0.3092 & 0.4287\\
            $Loss_f$ & median & 0.1882 & 0.3499 & 0.2333 & 0.2875 & 0.3368 & 0.3722 & 0.4694 & 0.514\\
                & mad & 0.0208 & 0.0352 & 0.0382 & 0.0267 & 0.0486 & 0.0579 & 0.0847 & 0.0786\\
            $Loss_{Bf}$ & median & 0.9577 & 1.9817 & 2.0844 & 2.6544 & 2.1681 & 2.4069 & 3.1446 & 3.4671\\
                & mad & 0.0902 & 0.188 & 0.287 & 0.2052 & 0.1734 & 0.2188 & 0.2869 & 0.3566\\
            $PrErr$ & median & 0.8923 & 1.9059 & 2.6905 & 3.1018 & 2.8433 & 3.1436 & 4.3324 & 4.9193\\
                & mad & 0.0902 & 0.1798 & 0.2362 & 0.3612 & 0.3206 & 0.395 & 0.354 & 0.5366\\
            \hline
        \end{tabular}
        %}
\end{table}

\begin{table}[h]
        \centering
        \caption{Simulation results: means and standard deviations of the four sample losses calculated for Thompson's factor scores over 100 runs.}
        \label{res_allT}
        %\resizebox{8cm}{3cm} {
        \begin{tabular}{|c|c|c|c|c|c|c|c|c|c|c|c|c|c|c|}
            \cline{1-10}
             & & \multicolumn{2}{|c|}{Setting 1} & \multicolumn{2}{|c|}{Setting 2} & \multicolumn{2}{|c|}{Setting 3} & \multicolumn{2}{|c|}{Setting 4}\\
            \cline{3-10}
            & & UNALCE & POET & UNALCE & POET & UNALCE & POET & UNALCE & POET\\
            \hline
            $Loss_B$ & mean & 2.8362 & 3.186 & 4.5014 & 4.701 & 3.5745 & 3.773 & 4.5492 & 4.8756\\
                & std & 0.1049 & 0.1586 & 0.1404 & 0.1829 & 0.2106 & 0.2564 & 0.4065 & 0.5361\\
            $Loss_f$ & mean & 0.1924 & 0.3566 & 0.2472 & 0.2916 & 0.3366 & 0.3848 & 0.4915 & 0.5305\\
                & std & 0.0265 & 0.0438 & 0.0798 & 0.0344 & 0.0631 & 0.073 & 0.1103 & 0.1003\\
            $Loss_{Bf}$ & mean & 0.9613 & 2.0299 & 2.1764 & 2.6577 & 2.1963 & 2.424 & 3.1759 & 3.5572\\
                & std & 0.1179 & 0.2435 & 0.5321 & 0.2565 & 0.2227 & 0.2749 & 0.3802 & 0.4871\\
            $PrErr$ & mean & 0.9064 & 1.921 & 2.674 & 3.2001 & 2.8525 & 3.1922 & 4.4129 & 5.0542\\
                & std & 0.1192 & 0.2277 & 0.2927 & 0.4532 & 0.3206 & 0.395 & 0.4614 & 0.7066\\
            \hline
        \end{tabular}
        %}
\end{table}

\begin{table}[h]
        \centering
        \caption{Simulation results: medians and median absolute deviations of the four sample losses calculated for Thompson's factor scores over 100 runs.}
        \label{res_allT2}
        %\resizebox{8cm}{3cm} {
        \begin{tabular}{|c|c|c|c|c|c|c|c|c|c|c|c|c|c|c|}
            \cline{1-10}
             & & \multicolumn{2}{|c|}{Setting 1} & \multicolumn{2}{|c|}{Setting 2} & \multicolumn{2}{|c|}{Setting 3} & \multicolumn{2}{|c|}{Setting 4}\\
            \cline{3-10}
            & & UNALCE & POET & UNALCE & POET & UNALCE & POET & UNALCE & POET\\
            \hline
            $Loss_B$ & median & 2.8471 & 3.1894 & 4.4875 & 4.6756 & 3.5586 & 3.7427 & 4.4653 & 4.7674\\
                & mad & 0.0854 & 0.1258 & 0.1122 & 0.1492 & 0.1704 & 0.2085 & 0.3084 & 0.4287\\
            $Loss_f$ & median & 0.1876 & 0.3499 & 0.2322 & 0.2875 & 0.3364 & 0.3722 & 0.4682 & 0.514\\
                & mad & 0.0208 & 0.0352 & 0.0383 & 0.0267 & 0.0486 & 0.0579 & 0.0845 & 0.0786\\
            $Loss_{Bf}$ & median & 0.9506 & 1.9817 & 2.0748 & 2.6544 & 2.1691 & 2.4069 & 3.1316 & 3.4671\\
                & mad & 0.0910 & 0.188 & 0.2826 & 0.2052 & 0.1737 & 0.2188 & 0.2863 & 0.3566\\
            $PrErr$ & median & 0.8923 & 1.9059 & 2.6905 & 3.1018 & 2.8433 & 3.1436 & 4.3324 & 4.9193\\
                & mad & 0.0902 & 0.1798 & 0.2362 & 0.3612 & 0.3206 & 0.395 & 0.354 & 0.5366\\
            \hline
        \end{tabular}
        %}
\end{table}

\begin{figure}[htb]
\centering
\makebox{
\includegraphics[width=0.6\textwidth]{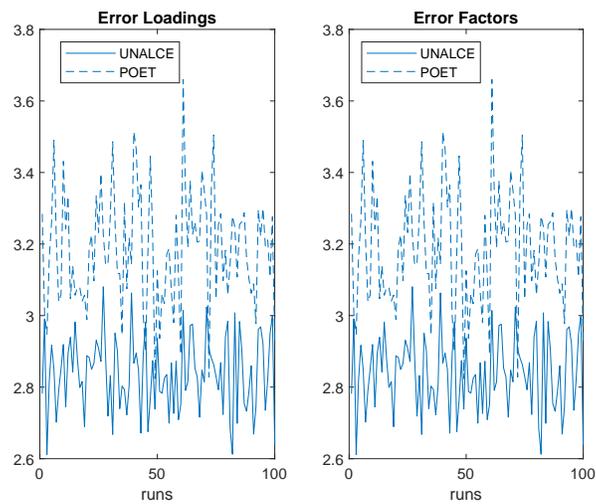}}
\caption{Bartlett's estimates: $Loss_B$, $Loss_f$ for Setting 1 over $100$ replicates.}
\label{fac_bart_1_1}
\end{figure}

\begin{figure}[htb]
\centering
\makebox{
\includegraphics[width=0.6\textwidth]{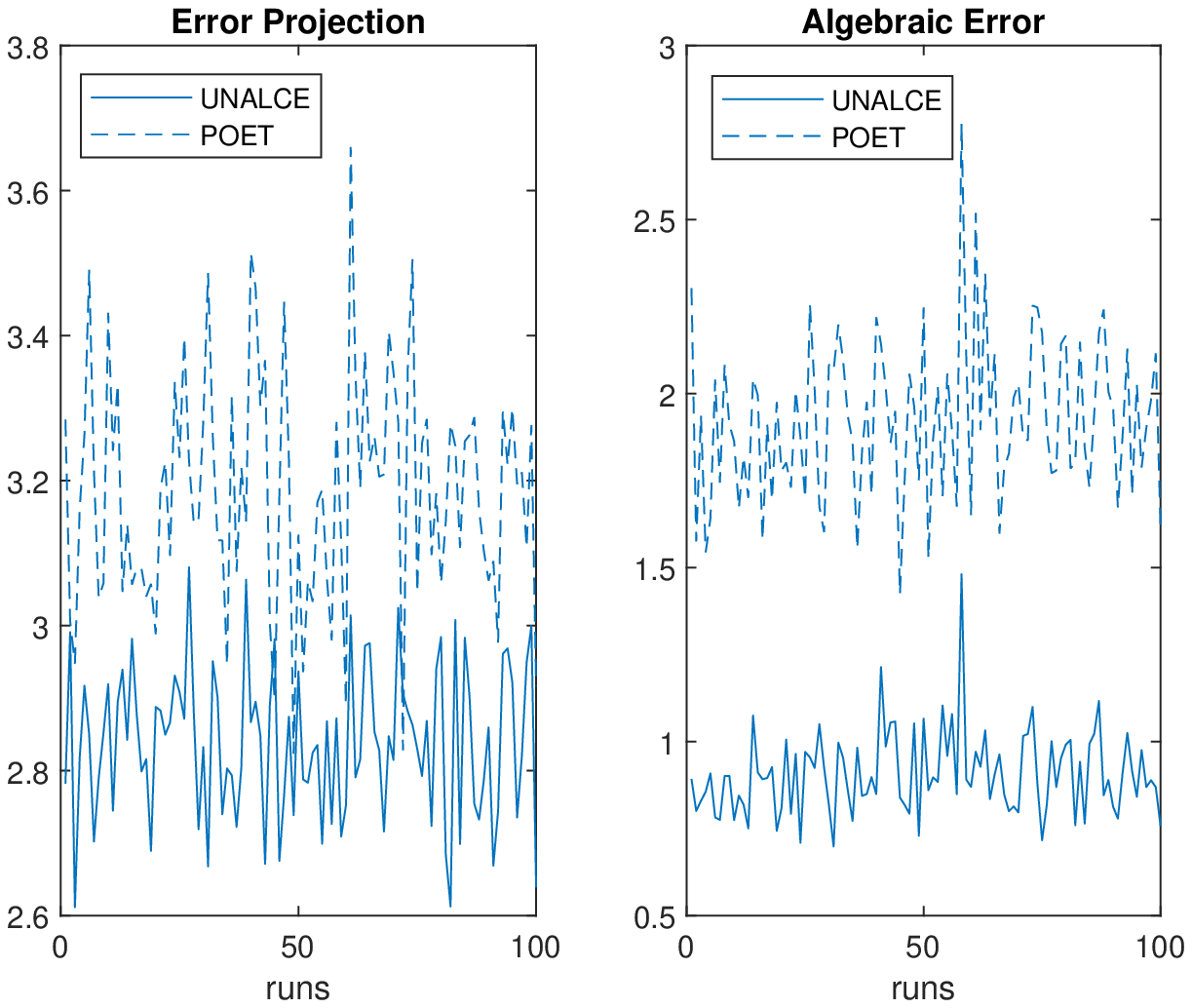}}
\caption{Bartlett's estimates: $Loss_{Bf}$ and $PrErr$ for Setting 1 over $100$ replicates.}
\label{fac_bart_1_2}
\end{figure}

\begin{figure}[htb]
\centering
\makebox{
\includegraphics[width=0.6\textwidth]{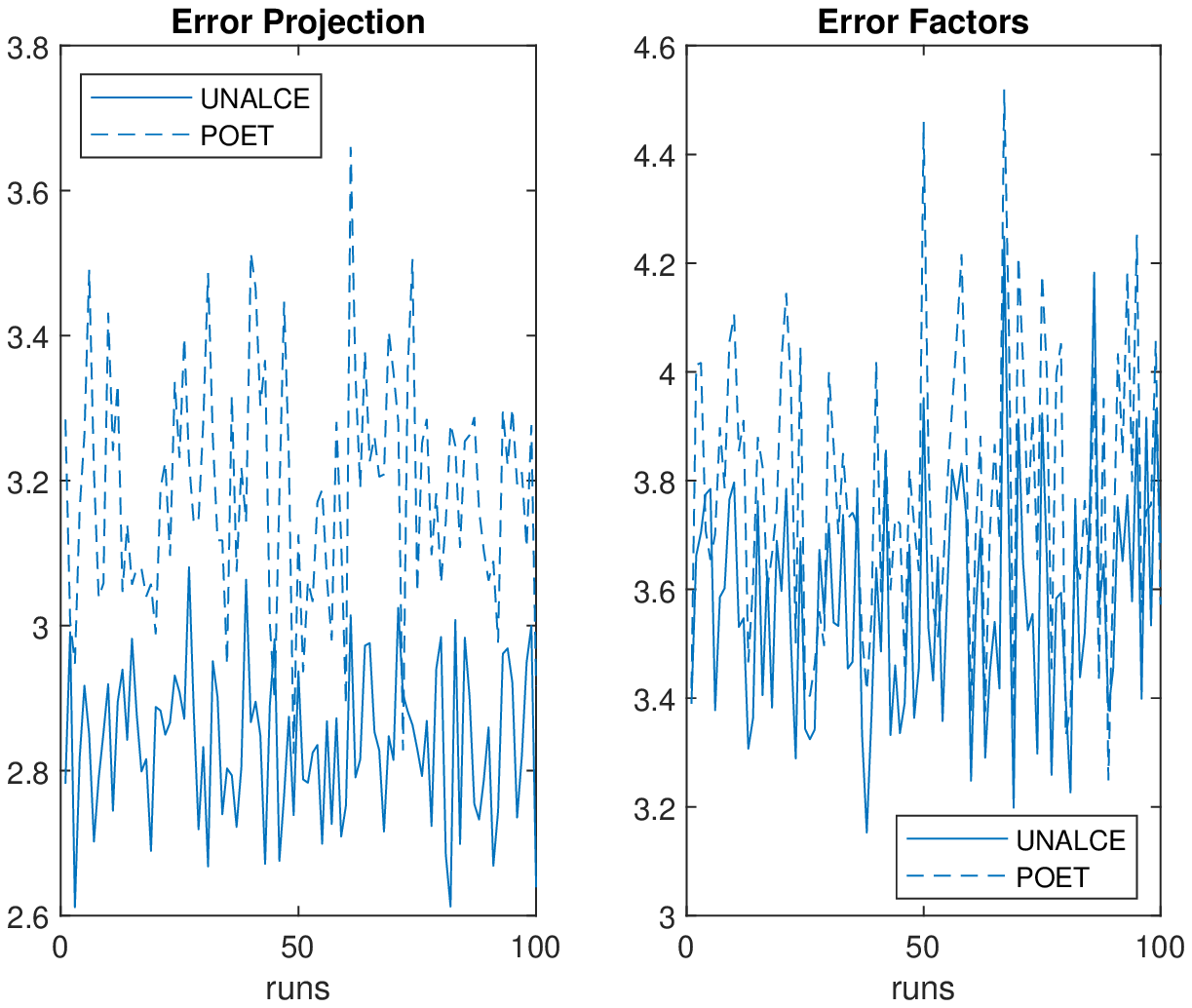}}
\caption{Bartlett's estimates: $Loss_B$, $Loss_f$ for Setting 4 over $100$ replicates.}
\label{fac_bart_4_1}
\end{figure}

\begin{figure}[htb]
\centering
\makebox{
\includegraphics[width=0.6\textwidth]{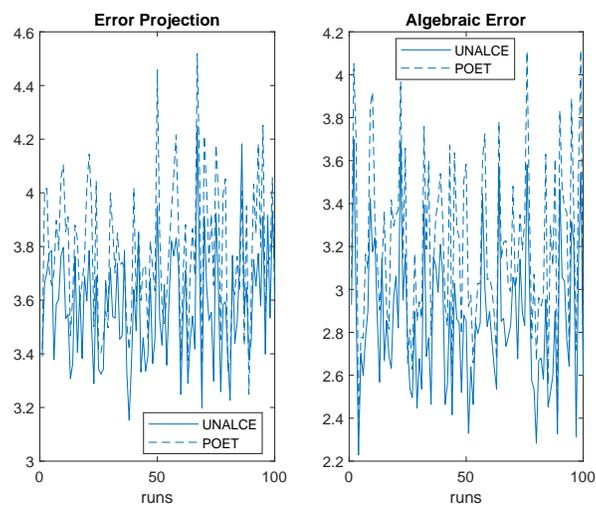}}
\caption{Bartlett's estimates: $Loss_{Bf}$ and $PrErr$ for Setting 4 over $100$ replicates.}
\label{fac_bart_4_2}
\end{figure}

In Table \ref{res_all}, we have reported means and standard deviations for the performance indicators
$Loss_B$, $Loss_f$, $Loss_{Bf}$ and $PrErr$, measured for Bartlett's factor scores
over 100 replicates for each setting.
We can observe that the means are smaller for UNALCE with respect to POET for each setting and indicator, while the variances tend to be larger, particularly for Settings 2 and 3.
This happens because Setting 2 has the most spiked eigenvalues and the smallest latent condition number, which leads to sporadic identifiability problems. As a proof of that, when we consider median and median absolute deviations, reported in Table \ref{res_all2},
UNALCE prevails over POET under all settings.
Anyway, we observe that the gain of UNALCE versus POET is far larger in Setting 1 and decreases progressively for Settings 2,3,4,
as those settings are increasingly consistent with POET assumptions.
This can also be appreciated in Figures \ref{fac_bart_1_1}, \ref{fac_bart_1_2}, \ref{fac_bart_4_1},
\ref{fac_bart_4_2}, which show $Loss_B$, $Loss_f$ and $Loss_{Bf}$, $PrErr$
for Settings 1 and 4 respectively.

%tables $5 \times (3+1)$
%$5 plots \times 4$!

\clearpage

\section{A real data example}\label{sec:real}

In this section, we apply the UNALCE methodology to a real financial dataset, already used
in \cite{fan2013large} to describe the performance of POET methodology.
The dataset contains $251$ annualized daily returns (year $2010$) of $p=50$ stocks, relative to five UK industry sectors: ''consumer goods-textiles and apparel clothing'', ''financial-credit services'', ''healthcare-hospitals'', ''services-restaurant'' and ''utilities-water utilities'',
with $10$ stocks from each sector.
%\vspace{0.2cm}

\begin{figure}[htb]
\centering
\makebox{\label{eig7}
\includegraphics[width=0.5\textwidth]{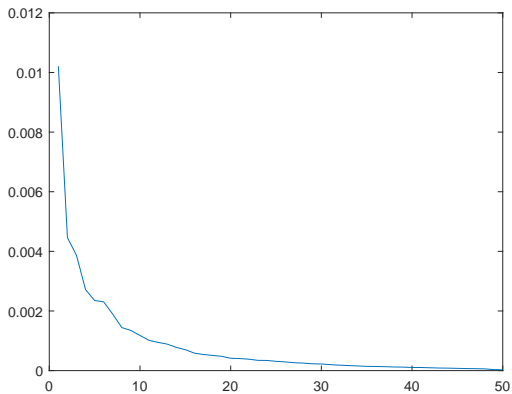}}
\caption{Sample eigenvalues: UK market data}
\end{figure}

In Figure \ref{eig7}, we report the eigenvalues of the $50-$dimensional sample covariance matrix.
Looking at the figure from a factor model perspective,
we can state that no more than $3$ latent factors should be considered.

In \cite{fan2013large}, the authors show the results of POET methodology with $r=2$, reporting that $25.8\%$ of recovered residual non-zeros are within blocks, and $6.7\%$ are off-blocks. In addition, all recovered non-zeros are positive within blocks, while only $60.3\%$ are positive off-blocks. The results are claimed to be similar for $r=1,2,3$.

In order to make a comparison, we have computed UNALCE estimates for a grid of $20 \times 20$ thresholds. The statistics of the optimal solutions, recovered via MC criterion (see \cite{farne2020large}), are reported in Table \ref{uk_res_1} (the optimal thresholds are $\widehat{\psi}=0.0007$ and $\widehat{\rho}=0.0004$).
We can note that the recovered rank is $1$, the UNALCE proportion of latent variance is very low (under $20\%$), and the recovered residual is diagonal.

Since UNALCE, differently from POET, recovers exactly the rank and the sparsity pattern, we have computed POET solutions (with hard thresholding) for $r=1$, selecting via $10-$fold cross-validation the optimal constant $\widehat{C}=1.10$ over a grid of $1000$ constants, linearly spaced from $0$ to $100$. The results are reported in Table \ref{uk_res_1}. We can note that the latent variance proportion is still very low ($23.29\%$), and the residual presents a relevant proportion of off-diagonal non-zeros ($19\%$) but a very small proportion of residual covariance ($0.89\%$). This means that the recovered non-zeros are irrelevant to explain the covariance structure. What is more, only $10.2\%$ of within-blocks elements are non-zeros, against the $21\%$ of off-blocks elements, and all recovered non-zeros are positive.

\begin{table}%\label{uk_res_1}
        \centering
        \caption{Covariance estimation results for UNALCE and POET (with $r=1$). $\widehat{r}$ is the latent rank, $\widehat{\theta}$ is the latent variance proportion, $\widehat{\rho}_{\widehat{S}}$ is the residual covariance proportion, $\widehat{\pi}_{nz}$ is the residual nonzero proportion, $\Vert \widehat{\Sigma}-\Sigma_n \Vert$ is the sample total loss.}
        \label{uk_res_1}
        %\resizebox{5cm}{2cm} {
        %\resizebox{8cm}{3cm} {
        \begin{tabular}{|c|c|c|}
            \hline
            & UNALCE & POET\\
            \hline
            $\widehat{r}$ & 1 & 1\\
            $\widehat{\theta}$ & 0.1930 & 0.2329 \\
            $\widehat{\rho}_{\widehat{S}}$  & 0 & 0.0089\\
            $\widehat{\pi}_{nz}$  & 0 & 0.1902 \\
            $\Vert \widehat{\Sigma}-\Sigma_n \Vert$ & \textbf{0.0013} & 0.0021\\
            \hline
        \end{tabular}
\end{table}

%The UNALCE method applied on the same data recovers a latent rank equal to one.
%In Table \ref{uk_res}, we compare the results of the two methods.
%%calculate factor model estimates. The results are stored .
%\begin{table} \label{uk_res}
%        \centering
%        \caption{Covariance estimation results of UNALCE and POET. $\widehat{r}$ is the latent rank, $\widehat{\theta}$ is the latent variance proportion, $\widehat{\rho}_{\widehat{S}}$ is the residual covariance proportion, $\widehat{\pi}_{nz}$ is the residual nonzero proportion, $\Vert \widehat{\Sigma}-\Sigma_n \Vert$ is the sample total loss.}
%        \label{uk_res}
%        %\resizebox{5cm}{2cm} {
%        \label{res_all2}
%        %\resizebox{8cm}{3cm} {
%        \begin{tabular}{|c|c|c|}
%            \hline
%            & UNALCE & POET\\% & UNALCE & POET & UNALCE & POET & UNALCE & POET\\
%            \hline
%            $\widehat{r}$ & 1 & 2\\
%            $\widehat{\theta}$ & 0.1889 & 0.3307 \\
%            $\widehat{\rho}_{\widehat{S}}$  & 0.0092 & 0.0648 \\
%            $\widehat{\pi}_{nz}$  & 0.0122 & 0.1069\\
%            $\Vert \widehat{\Sigma}-\Sigma_n \Vert$ & \textbf{0.0023} & 0.0028\\
%            \hline
%        \end{tabular}
%\end{table}
%%%go!

Since we know that POET does not offer any algebraic guarantee on the recovered sparsity pattern, and UNALCE approximates quite better the sample covariance matrix (the sample total loss is $0.0013$ against $0.0021$ for POET),
we cannot claim so easily the presence of a residual cluster-wise structure.
On the contrary, it is more likely that we have one weak latent factor (consistently recovered by UNALCE)
which entirely explains the covariance structure.
The result reported in \cite{fan2013large}
could be explained by the use of principal component analysis with $p=50$,
which could be not large enough to ensure that the estimated residuals are no longer correlated across variables.

Finally, we calculate three quantities to estimate the variability of the estimated loadings, factor scores and factor projections:
\begin{itemize}
\item $var_B=\sum_{j=1}^p \Vert \widehat{b}_{j}-\tilde{b}\Vert$,
where $\tilde{b}$ is the mean estimated factor loading;
\item $var_f=\sum_{k=1}^n \Vert \widehat{f}_{k}\Vert$ for factor scores ($\sum_{k=1}^{n}{\widehat{f}_{k}}$ is zero by construction);
\item $var_{Bf}=\sum_{k=1}^n \Vert \widehat{B} \widehat{f}_{k}- \widetilde{Bf}\Vert$, where $\widetilde{Bf}$ is the mean factor projection across the observations.
\end{itemize}

These computations are reported for Bartlett's and Thompson's UNALCE and POET factor model estimates
in Table \ref{uk_losses}.
%observe that the quantity $var_B=\sum_{k=1}^n \Vert \widehat{f}_{k}\Vert$ is smaller for UNALCE than for POET. The same holds for the variability of estimated loadings, expressed as $var_f=\sum_{j=1}^p \Vert \widehat{b}_{j}-\bar{\widehat{b}}\Vert$, where $\bar{\widehat{b}}$ is the mean factor loading vector, and the variability of the projections, estimated by $var_{Bf}=\sum_{k=1}^n \Vert \widehat{B} \widehat{f}_{k}- \bar{\widehat{B} \widehat{f}}\Vert$, where $\bar{\widehat{B} \widehat{f}}$ is the mean projection.
We observe that $var_B$ and $var_{Bf}$ are better for UNALCE, as we could expect from Corollary \ref{Coroll_B} and Theorems \ref{bartlett_opt} and \ref{thompson_opt}, while $var_f$ presents smaller values for POET. This is due to the particular structure of the POET residual component, which presents only positive elements. However, we must note that such structure cannot be trusted, as UNALCE recovers a diagonal residual and possesses the algebraic consistency property.

%For UNALCE, we have that $29.39\%$ of the covariance is thrown away. %, while $78.15\%$ lies in the low rank component and $0.72\%$ in the sparse component.
%The covariance structure appears so complex that a relevant proportion
%of \textbf{residual covariance} is present.
%For POET, the proportion of \textbf{latent variance} is much \textbf{larger}
%while we find \textbf{much less} residual covariance.
%The UNALCE estimation process approximates \textbf{better} than POET the sample covariance matrix.

\begin{table}[h]
        \centering
        \caption{Estimated variabilities for UNALCE and POET ($r=1$). Bartlett's and Thompson's factor loadings, scores and projections.}
        \label{uk_losses}
        %\resizebox{8cm}{3cm} {
        \begin{tabular}{|c|c|c|c|c|c|c|c|c|c|c|c|c|c|c|}
            %\cline{1-10}
            % & & \multicolumn{2}{|c|}{Setting 1} & \multicolumn{2}{|c|}{Setting 2} & \multicolumn{2}{|c|}{Setting 3} & \multicolumn{2}{|c|}{Setting 4}\\
            \cline{1-4}
            Metric & Method & UNALCE & POET \\
            \hline
            $var_B$ &  & 0.1990 & 0.2367 \\
            $var_f$ & Bartlett & 197.53 & 194.25 \\
                & Thompson & 189.74 & 174.65 \\
            $var_{Bf}$ & Bartlett & 18.17 & 19.63 \\
                & Thompson & 17.46 & 17.65\\
            \hline
        \end{tabular}
        %}
\end{table}

\section{Conclusions}\label{sec:concl}

In this paper, we propose a new method to estimate an approximate factor model
with a sparse residual in high dimensions.
In particular, we elaborate over the results of \cite{farne2020large} to prove that the ordinary least squares (OLS) estimates of factor loadings and scores based on UNALCE (UNshrunk ALgebraic Covariance Estimator) are asymptotically consistent, as well as the same estimates based on POET \citep{fan2013large}.
Consistency holds in Euclidean norm under the assumption of intermediate spikiness of latent eigenvalues
and element-wise sparsity of the residual component,
while UNALCE provides the exact recovery of the latent rank and the residual sparsity pattern.
A lower bound is imposed  on the smallest latent eigenvalue and the smallest absolute nonzero residual
element to ensure identifiability.

Moving from the eigenvalue dispersion lemma of \cite{ledoit2004well}, we then prove that Bartlett's and Thompson's factor scores show the tightest possible error bound in Euclidean norm given the finite sample,
within the class of estimate pairs with exact low rank and sparsity pattern.
The proofs require advanced techniques of matrix algebra. In addition, it is proved that the projection of the low rank error matrix onto the orthogonal complement of the low rank space has the minimum possible Euclidean norm given the finite sample.
Moreover, Bartlett's and Thompson's scores converge to the OLS ones, thus being also asymptotically consistent.

In the end, we prove in an ad hoc simulation study the validity of our optimality results, showing that UNALCE-based factor scores work particularly well with respect to POET-based ones if the latent eigenvalues are not so spiked and the residual is very sparse with prominent non-zeros in absolute value.
A real financial data example further supports the optimality properties of the UNALCE approach compared to the POET one.
%large smallest latent eigenvalue
%large minimum absolute nonzero element

%\newpage

\appendix

\section{Assumptions and key results of \cite{farne2020large}}\label{ass:fm20}

\subsection{Assumptions}

\begin{Ass}\label{eigenvalues}
All the eigenvalues of the $r\times r$ matrix $p^{-\alpha}\mf{B}^\top \mf{B}$
are bounded away from $0$ for all $p$ and $\alpha \in [0,1]$.
\end{Ass}

%In order to guarantee identifiability using (\ref{ours}), the following assumption is crucial.
\begin{Ass}\label{alg}
There exist $k_L,k_S>0$, $\delta \in [0,0.5]$, such that $\xi(T(\mf{L}^{*}))=\sqrt{{r}/{(k_L^2 p^{2\delta})}}$, $\mu(\Omega(\mf{S}^{*}))=k_S p^{\delta}$, ${k_S}/{k_L}\leq {1}/{54}$ with $\delta < \alpha$.
\end{Ass}

\begin{Ass}\label{tails}
There exist $r_1,r_2>0$ and $b_1,b_2>0$ such that, for any $t>0$, $k\leq n$, $i \leq r$, $j\leq p$:
\begin{eqnarray}
\Pr(\vert\vf{f}_{ik}\vert>s) \leq \exp({-b_1/t}), \qquad \Pr(\vert\vf{\epsilon}_{jk}\vert>s) \leq \exp({-b_2/t}). \nonumber
\end{eqnarray}
\end{Ass}
%$\log(p)=o[(3 b_1+1.5 b_2+1)^{-1}]$.

\begin{Ass}\label{sparsity}
There exist $c_1,c_2,c_3,\delta_2,\delta_2'>0$, $\delta' \in [0,\delta + 0.5]$
such that $\lambda(\mf{S}^*)_{min} > c_1$, $\min_{i,j \leq p}\mathrm{var}(\vf{\epsilon}_{ik}\vf{\epsilon}_{jk})>c_2$
for any $k \leq n$, $i,j \leq p$, $s^{*}_{ii} \leq c_3$ for any $i\leq p$,% p^{\delta}
$s'=\max_{i \leq p} \sum_{j \leq p} \mathbbm{1}(\mf{S}^*_{ij}=0)\leq \delta_2 p^{\delta}$ with $\delta_2\geq k_S$ and
$\Vert \mf{S}\Vert_{1,v}=\max_{i \leq p} \sum_{j \leq p} \vert\mf{S}^*_{ij}\vert \leq \delta_2' p^{\delta'}$.
%and $\delta_4>0$.
\end{Ass}

\begin{Ass}\label{pr}
There exist $\delta_3,\delta_4>0$ such that $r = \delta_3 \ln{p}$ and $n\geq \delta_4 p^{1.5\delta}$.
\end{Ass}

\begin{Ass}\label{sc}
$2\delta \leq \alpha \leq 2\delta+\delta'$ and $\delta_5<\frac{C_3 \delta}{k_L\delta_4}<\delta'$ with $\delta_5>0$.
\end{Ass}

Assumption \ref{eigenvalues} prescribes
that the latent eigenvalues are spiked
in the sense of Yu and Samworth (\cite{fan2013large}, p. 656),
thus allowing for intermediately pervasive latent factors as $p$ diverges.
Assumption \ref{alg} ensures the identifiability of $L^{*}$ and $S^{*}$.
Note that the condition $\delta < \alpha$ imposes a gap between the magnitude of latent eigenvalues and
the residual sparsity degree, i.e. the number of residual non-zeros.
Assumption \ref{tails} bounds the tails of factors and residuals. Assumption \ref{sparsity}
controls for the sparsity degree of the residual component.
Assumption \ref{pr} prescribes the order of the latent rank
and a necessary lower bound for the sample size.
Assumption \ref{sc} ensures the sparsity pattern recovery.
We refer to \cite{farne2020large} for more details.

\subsection{Key results}

\begin{Thm}\label{thmMinetop}
Let $\Omega=\Omega({S}^{*})$ and
$T=T({L}^{*})$. Suppose that Assumptions \ref{eigenvalues}-\ref{sc} and Assumption \ref{lowerbounds} hold.
%$m_p\ll C(p^{\alpha})$ for $q=0$,
%\textbf{IF}:\\
%1) exactly
Define
$$
\psi=\frac{1}{\xi(T)}\frac{p^{\alpha}}{\sqrt{n}}
%\lambda=
%\left\{
%\begin{array}{rl}
% & 0.5 \leq \alpha \leq 1,\\
%\frac{1}{\xi(T)}\frac{p^{\alpha}}{\sqrt{n}} & 0\leq \alpha < 0.5,
%\end{array}
%\right.
$$
with $\rho=\gamma \psi$, where $\gamma \in [9\xi(T),1/(6\mu(\Omega))]$.
%\textbf{Identification} conditions:\\
%2)$\mu(\Omega)\xi(T)\leq \frac{1}{54}$;\\
%3)\\
%%1) and 2) are needed for \textbf{identification}.\\ %(impose a bound on the curvature of T).\\
%4)
%5)
%6)
%$n\geq p$ (needed for ${\Sigma}_n$).\\ %%possibly different!
Then, with probability greater than
$1-C_4p^{-C_5}$, %\\\vspace{0.1cm}
the pair $(\widehat{{L}},\widehat{{S}})$ minimizing (\ref{ours})
recovers the rank of ${L}^{*}$ ($rank(\widehat{{L}})=rank({L}^{*})$)
and the sparsity pattern of $S^{*}$
($sign(\widehat{{S}})=sign({S}^{*})$).
Moreover, the matrix losses for each component are bounded as follows:
$$\Vert\widehat{{L}}-{L}^{*}\Vert_2\leq
C\psi, \qquad \Vert\widehat{{S}}-{S}^{*}\Vert_{\infty}\leq C\rho.$$
\end{Thm}

Once defined $\phi_S=C s' \xi(T) \psi$ and $\phi=C(s'\xi(T)+1)\psi$,
where $s'$ is the maximum number of non-zeros per row/column in ${S}^{*}$,
we can state the following Corollary.
\begin{Coroll}\label{losses}
Under the assumptions of Theorem \ref{thmMinetop}, it holds
$||\widehat{{S}}_{ALCE}-{S}^{*}||_2~\leq \phi_S $, $||\widehat{{\Sigma}}_{ALCE}-{\Sigma}^{*}||_{2} \leq \phi$, $||\widehat{{S}}_{ALCE}^{-1}-{S}^{*-1}||_2 \leq \phi_S$,
and $||\widehat{{\Sigma}}_{ALCE}^{-1}-{\Sigma}^{^*-1}||_{2} \leq \psi$.
In addition, $\widehat{{S}}_{ALCE}$, $\widehat{{\Sigma}}_{ALCE}$,  $\widehat{{S}}_{ALCE}^{-1}$,  $\widehat{{\Sigma}}_{ALCE}^{-1}$ are positive definite if and only if $\lambda_{p}({{S}^{*}})>\phi_{{S}}$, $\lambda_{p}({{\Sigma}^{*}})>\phi$, $\lambda_{p}({{S}^{*}}) \geq 2 \phi_{{S}}$,  $\lambda_{p}({{\Sigma}^{*}}) \geq 2 \phi$, respectively.
\end{Coroll}

We refer to \cite{farne2020large} for the proofs.

\section{Proofs}\label{proofs}

\subsection{Proof of Theorem \ref{facPOET}}\label{facFan}

%Suppose that we decompose ${E}_n={\Sigma}_n-{\Sigma}^*$ in its four components:
%$${E}_n={\Sigma}_n-{\Sigma}^*={D}_1+{D}_2+{D}_3+{D}_4,$$ where:
%$${D}_1=\left(n^{-1} {B} \sum_{k=1}^n {f}_k {f}_k'-{I}_r\right){B}',$$
%$${D}_2=n^{-1} \left(\sum_{k=1}^n {\epsilon}_k {\epsilon}_k'-{S}^*\right),$$
%$${D}_3= n^{-1}{B}\sum_{k=1}^n {f}_k {\epsilon}_k',$$
%$${D}_4={D}_3',$$
%where ${{f}}_k$ and ${{\epsilon}}_k$ are respectively the vectors of factor scores and residuals for each observation.
%%Asterisks are omitted to avoid cluttered notation.

Recalling that ${\Sigma}_n=(n-1)^{-1}\sum_{k=1}^n {x}_k {x}_k'$ and ${x}_k={B}{f}_k+{\epsilon}_k$,
where ${{f}}_k$ and ${{\epsilon}}_k$, $k=1,\ldots,n$, are respectively the vectors of factor scores and residuals for each observation, we can decompose the error matrix ${{E}}_n={\Sigma}_n-{\Sigma}^*$ in four components as follows
(cf. \cite{fan2013large}):
$${{E}}_n={\Sigma}_n-{\Sigma}^*=\widehat{{D}}_1+\widehat{{D}}_2+\widehat{{D}}_3+\widehat{{D}}_4,$$ where:
\begin{equation}
\widehat{{D}}_1=\left(n^{-1} {B} \sum_{k=1}^n {f}_k {f}_k'-{I}_r\right){B}',
\widehat{{D}}_2=n^{-1} \sum_{k=1}^n \left({\epsilon}_k {\epsilon}_k^\top-{S}^*\right),
\widehat{{D}}_3= n^{-1}{B}\sum_{k=1}^n {f}_k {\epsilon}_k',
\widehat{{D}}_4=\widehat{{D}}_3'.\nonumber
\end{equation}
We thus recall from \cite{farne2020large} the following result.
\begin{Lemma}\label{Lemma4}
Under Assumptions \ref{eigenvalues}, \ref{tails}, \ref{sparsity}, \ref{pr},
there exists a positive constant $C$ such that
\begin{eqnarray}
||{D}_1||_2 &\leq& C \left(p^{\alpha}\sqrt{\frac{1}{n}}\right); \nonumber\\
||{D}_2||_{2} &\leq& C \left(p^{\delta}\sqrt{\frac{\log{p}}{n}}\right); \nonumber \\
||{D}_3||_{2} &\leq& C\left(p^{\frac{\alpha}{2}+\frac{\delta}{2}}\sqrt{\frac{\log{p}}{n}}\right). \label{prodlemma}
\end{eqnarray}
\end{Lemma}

As a consequence, we can recall the following fundamental Lemma proved in \cite{farne2020large}. %follows
%, and recalling Lemma 5 in \cite{fan2013large}, it is proved
\begin{Lemma}\label{samplego}
Let $\widehat{\lambda}_r$ be the $r-$th largest eigenvalue of ${{\Sigma}}_n$. Under the assumptions of Lemma \ref{Lemma4},
if $\delta < \alpha$, then $\widehat{\lambda}_r>C_1p^{\alpha}$ with probability approaching $1$ for some $C_1>0$.
\end{Lemma}
%$E_n=\Sigma_n-\Sigma^{*}=D_1+D_2+D_3+D_4$.
%From \cite{farne2018finite}, we inherit the following Lemma.
%\begin{Lemma}\label{samplego}
%Let $\widehat{\lambda}_r$ be the $r-$th largest eigenvalue of $\widehat{{\Sigma}}_n$.
%Then $\widehat{\lambda}_r>C_1p^{\alpha}$ with probability approaching $1$ for some $C_1>0$.
%\end{Lemma}
%
%
%$$||D_1|| \leq O\left(\frac{p^{\alpha}}{\sqrt{n}}\right)$$
%$$||{\widehat{S}_n-S}||_{\infty}\leq O\left(p^{\delta-1}\sqrt{\frac{\log{p}}{n}}\right)$$
%$$||D_2||=O\left(\frac{p^{\delta}}{\sqrt{n}}\right)$$
%$$||D_3||\leq O\left(\frac{p^{\frac{\alpha}{2}+\frac{\delta}{2}}}{\sqrt{n}}\right)$$
%The theorem states that
%If \textbf{Ass.} \ref{4'} holds in place of \textbf{Ass.} \ref{4_1} and the \textbf{true} $r$ is known, it is possible to prove that
%$$\max_{j\leq p}||{\widehat{b}_j-H b_j}||=O_p\left(\omega_n\right)$$
%with $\omega_n=p^{\alpha}\sqrt{\frac{\log p}{n}}+\frac{1}{p^{\alpha/2}}$ and
%$$\max_{i\leq n}||{\widehat{f}_i-H f_i}||=O_p\left(\frac{p^{1}{n^{\frac{1}{2}}}+\frac{n^{1/4}}{p^{\alpha/2}}\right)$$
%$$\max_{j\leq p,i\leq n}||{\widehat{b}_j'\widehat{f}_i-b_j' f_i}||=O_p\left(p^{\alpha\sqrt{\frac{\log p}{n}}\log{n}^{1/b_2}+\frac{n^{1/4}}{p^{\alpha/2}}\right).$$
%\pause
%What about $\widehat{b}_{j,UNALCE}$, $\widehat{f}_{i,UNALCE}$, $\widehat{b}_{j,UNALCE}\widehat{f}_{i,UNALCE}$?

Then, we proceed as follows. According to \cite{bai2003inferential}, setting
$\widehat{F}=\widehat{F}_{OLS,1}$, we can write, for each $t=1,\ldots,n$,
$$\widehat{f}_t-H f_t=\left(\frac{\widehat{\Lambda}_r}{p}\right)^{-1}\left\{ \frac{1}{n} \sum_{s=1}^n \widehat{f}_s \frac{E(\epsilon_s'\epsilon_t)}{p}+\frac{1}{n} \sum_{s=1}^n \widehat{f}_s\varsigma_{st}+\frac{1}{n} \sum_{s=1}^n \widehat{f}_s\eta_{st}+\frac{1}{n} \sum_{s=1}^n \widehat{f}_s\xi_{st}\right\},$$
where $\xi_{st}=\frac{(f_t'\sum_{i=1}^p b_i \epsilon_{is})}{p}$, $\eta_{st}= \frac{(f_s'\sum_{i=1}^p b_i \epsilon_{it})}{p}$,  $\varsigma_{st}= \frac{\epsilon_s'\epsilon_t}{p}- \frac{E(\epsilon_s'\epsilon_t)}{p}$, and
$\widehat{\Lambda}_r$ is the diagonal matrix containing the top $r$ eigenvalues of $\Sigma_n$ in decreasing order. The eigenvalues of $\left(\frac{\widehat{\Lambda}_r}{p}\right)^{-1}$ scale to $O(p^{1-\alpha})$ as $p\rightarrow \infty$.
%Therefore, if $r$ is known, POET methodology is consistent even if $\alpha < 1$.

%However, if we replace
%$\widehat{\Lambda}_r$ by $\widehat{\Lambda}_{UNALCE}$ we proved in Section \ref{sec:eig_conc} that the estimated eigenvalues
%are the most concentrated possible around their true mean.
%under the prescribed assumptions for $UNALCE$.
%The most bounded possible in a finite sample!

%At the same time, the lower bound of $\left(\frac{\widehat{\Lambda}_r}{p}\right)^{-1}$
%and $\left(\frac{\widehat{\Lambda}_{ALCE}}{p}\right)^{-1}$ is the same.
%In fact, defined $\Sigma_r$ as the covariance matrix formed by the top $r$ principal components of $\Sigma_n$, it holds
%\begin{eqnarray}
%\lambda_r(\widehat{{L}}_{ALCE})={\lambda}_r(\widehat{{L}}_{ALCE}-{\Sigma}_r+{\Sigma}_r)\geq \nonumber\\
%\geq{\lambda}_{p-r+1}(\widehat{{L}}_{ALCE}-{\Sigma}_r)+\lambda_r({\Sigma}_r).\nonumber
%\end{eqnarray}
%Since ${\lambda}_{p-r+1}(\widehat{{L}}_{ALCE}-\Sigma_r)$ is asymptotically zero as $p^{\alpha}\log(p) \gg n \geq p^{\frac{3}{2}\delta}$,
%it follows that asymptotically $\big \vert \big \vert \left(\frac{\widehat{\Lambda}_r}{p}\right)^{-1}\big \vert \big \vert =\big \vert \big \vert \left(\frac{\widehat{\Lambda}_{ALCE}}{p}\right)^{-1} \big \vert \big \vert =O(p^{1-\alpha})$.

Relying on Assumption \ref{4'}, we can prove the following lemmas by simply applying the corresponding proofs in \cite{fan2013large}, where Assumption \ref{4'} is imposed with $\alpha=1$.
\begin{Lemma}\label{Lemma8}
For all $i\leq r$:
\begin{eqnarray}
\frac{1}{n} \sum_{t=1}^n \left(\frac{1}{n p} \sum_{s=1}^n \widehat{f}_{is} E(\epsilon_s'\epsilon_t)\right)^2&=&O(n^{-1});\nonumber\\
\frac{1}{n} \sum_{t=1}^n \frac{1}{n} \left(\sum_{s=1}^n \widehat{f}_{is} \varsigma_{st}\right)^2&=&O(p^{-\alpha});\nonumber\\
\frac{1}{n} \sum_{t=1}^n \frac{1}{n} \left(\sum_{s=1}^n \widehat{f}_{is} \eta_{st}\right)^2&=&O(p^{-\alpha});\nonumber\\
\frac{1}{n} \sum_{t=1}^n \frac{1}{n} \left(\sum_{s=1}^n \widehat{f}_{is} \xi_{st}\right)^2&=&O(p^{-\alpha}).\nonumber
\end{eqnarray}
\end{Lemma}

%Then, relying again on Assumption \ref{4'}, we can state the following Lemma.
\begin{Lemma}\label{Lemma9}
\begin{eqnarray}
\max_{t \leq n} \bigg \vert \bigg \vert\frac{1}{np} \sum_{s=1}^n \widehat{f}_s E(\epsilon_s'\epsilon_t)\bigg \vert \bigg \vert=O\left(\sqrt{\frac{1}{n}}\right);\nonumber\\
\max_{t \leq n} \bigg \vert \bigg \vert\frac{1}{n} \sum_{s=1}^n \widehat{f}_s \varsigma_{st}\bigg \vert \bigg \vert=O\left(\frac{n^{1/4}}{p^{\alpha/2}}\right);\nonumber\\
\max_{t \leq n} \bigg \vert \bigg \vert\frac{1}{n} \sum_{s=1}^n \widehat{f}_s \eta_{st}\bigg \vert \bigg \vert=O\left(\frac{n^{1/4}}{p^{\alpha/2}}\right);\nonumber\\
\max_{t \leq n} \bigg \vert \bigg \vert\frac{1}{n} \sum_{s=1}^n \widehat{f}_s \xi_{st}\bigg \vert \bigg \vert=O\left(\frac{n^{1/4}}{p^{\alpha/2}}\right).\nonumber
\end{eqnarray}
\end{Lemma}

We recall that $H=\frac{1}{n}(\widehat{\Lambda}_r)^{-1}\widehat{F}'FB'B$.
Since the eigenvalues of $\left(\frac{\widehat{\Lambda}_r}{p}\right)^{-1}$ scale to $O(p^{1-\alpha})$ instead of $O(p)$,
applying Assumption \ref{4'} and following \cite{fan2013large} we obtain

%\begin{Lemma}\label{Lemma9}
%\begin{eqnarray}
%\max_{i\leq r}  \frac{1}{n} \sum_{k=1}^n (\widehat{f}_k-Hf_k)_i^2=O(p^{1-\alpha})O\left(\frac{1}{n}+\frac{1}{d}\right)\nonumber\\
%\frac{1}{n} \sum_{k=1}^n ||\widehat{f}_k-Hf_k||=O(p^{1-\alpha})O\left(\frac{1}{n}+\frac{1}{d}\right)\nonumber\\
%\max_{k \leq n} \sum_{k=1}^n||\widehat{f_k}-H f_k||=O(p^{1-\alpha})O\left(\sqrt{\frac{1}{n}}+\frac{n^{1/4}}{p^{\alpha/2}}\right)\label{th2}
%\end{eqnarray}
%\end{Lemma}
%Note that (\ref{th2}) is the second thesis of Theorem \ref{facPOET}.

\begin{Lemma}\label{Lemma10}
\begin{eqnarray}
\max_{i\leq r}  \frac{1}{n} \sum_{k=1}^n (\widehat{f}_k-Hf_k)_i^2=O\left(\frac{p^{1-\alpha}}{n}+\frac{p^{1-\alpha}}{p^{\alpha}}\right);\nonumber\\
\frac{1}{n} \sum_{k=1}^n ||\widehat{f}_k-Hf_k||=O\left(\frac{p^{1-\alpha}}{n}+\frac{p^{1-\alpha}}{p^{\alpha}}\right);\nonumber\\
\max_{k \leq n} \sum_{k=1}^n||\widehat{f_k}-H f_k||=O\left({\frac{p^{1-\alpha}}{\sqrt{n}}}+\frac{p^{1-\alpha}n^{1/4}}{p^{\alpha/2}}.\right)\label{th2}
\end{eqnarray}
\end{Lemma}
Note that (\ref{th2}), which bounds the uniform rate of $\widehat{f}_k-H f_k$ over $k=1,\ldots,n$,
proves the second thesis of Theorem \ref{facPOET}.

At the same time, applying Assumption \ref{4'} and following \cite{fan2013large}, we can claim
\begin{Lemma}\label{Lemma11}
\begin{eqnarray}
HH'=I_{\widehat{r}}+O\left(\frac{p^{(1-\alpha)/2}}{\sqrt{n}}+\frac{p^{(1-\alpha)/2}}{p^{\alpha/2}}\right);\nonumber\\
H'H=I_{{r}}+O\left(\frac{p^{(1-\alpha)/2}}{\sqrt{n}}+\frac{p^{(1-\alpha)/2}}{p^{\alpha/2}}\right).\nonumber
\end{eqnarray}
\end{Lemma}

We can now prove the first thesis of Theorem \ref{facPOET},
about the uniform rate of $\widehat{b}_j-H b_j$ over $j=1,\ldots,p$.
Following \cite{fan2013large}, we know that $\widehat{b}_j-H b_j$ can be decomposed in three terms:
\begin{equation}
\widehat{b}_j-H b_j=\frac{1}{n} \sum_{k=1}^n H f_k \epsilon_{k,i} + \sum_{k=1}^n x_{k,i} (\widehat{f}_k -H f_k) + H (\widehat{f}_k\widehat{f}_k'-I_r)b_j=I+II+III.
\end{equation}
The claim (\ref{prodlemma}) in Lemma \ref{Lemma4} allows us to prove that $I$
is $O\left(p^{\delta/2}\sqrt{\frac{\log(p)}{n}}\right)$.
From Lemma \ref{Lemma10} and the fact $||H||=O(1)$, it follows that the $II$ is $O\left(\frac{p^{\frac{1-\alpha}{2}}}{n^{\frac{1}{2}}}+\frac{p^{\frac{1-\alpha}{2}}}{p^{\frac{\alpha}{2}}}\right)$.
Lemma \ref{Lemma10} and Assumption \ref{4'} imply that $III$ is $O\left(\frac{1}{\sqrt{n}}\right)$.
Therefore, we can prove that $$\max_{j\leq p}||{\widehat{b}_j-H b_j}||=O_p\left(\omega_n\right)$$
where $\omega_n=p^{\delta/2}\sqrt{\frac{\log p}{n}}+\frac{p^{(1-\alpha)/2}}{p^{\alpha/2}}$.

%$\omega_n=p^{1-\alpha}\sqrt{\frac{\log(p)}{n}}+\frac{1}{p^{\alpha/2}}$
%if $\alpha\leq 2/3+\frac{4}{3}\delta$
%$\omega_n=p^{\alpha/2-2\delta}\sqrt{\frac{\log(p)}{n}}+\frac{1}{p^{\alpha/2}}$
%if $2/3+\frac{4}{3}\delta<\alpha \leq 1$.

%Lemma 11
%$HH'=I_{\widehat{r}}+O(p^{1-\alpha})O(1/sqrt{n}+1/{(p^{\alpha/2})})$
%Note that it is still $||H||=O(1)$.
%$H'H=I_{{r}}+O(p^{1-\alpha})O(1/sqrt{n}+1/{(p^{\alpha/2})})$
%
%$\max_{j \leq r} ||\widehat{b}_j-Hb_j||=O((p^\alpha/2)\sqrt{\frac{\log(p)}{n}}+\sqrt(O(p^{1-\alpha})O(\frac{1}{n}+1/p^{\alpha/2}))$
%$\sqrt(O(p^{1-\alpha})O(\frac{1}{n}+1/p^{\alpha/2}))=p^{\alpha/2-\frac{1}{2}}/{\sqrt{(n}}+p^{\alpha/4-\frac{1}{2}}=\frac{p^{\alpha/2}}{\sqrt{n}}+\frac{p^{\alpha/4}}{\sqrt{p}}$
%$=O({p^{\alpha/2}})O({\sqrt{\frac{\log(p)}{n}}})}$.
%
%Finally, we can show that
%$\max_{j\leq p, i\leq r}||\widehat{b}_j' \widehat{f}_i-b_j'f_i||=O(\frac{p^{1-\alpha} n^{1/4}{p^{\alpha/2}}}+\log(n)^{1/b_2}p^{\alpha/2}\sqrt{\frac{\log{p}{n}}})$, if $\alpha\leq 2/3$, and
%$\max_{j\leq p, i\leq r}||\widehat{b}_j' \widehat{f}_i-b_j'f_i||=O(\frac{p^{1-\alpha} n^{1/4}{p^{\alpha/2}}}+\log(n)^{1/b_2}p^{1-\alpha}\sqrt{\frac{\log{p}{n}}})$, if $2/3<\alpha \leq 1$.

%$\max_{j \leq r} ||\widehat{b}_j-Hb_j||=$
%$=O(p^{\alpha/2-2\delta})\sqrt{\frac{\log(p)}{n}}+O(\frac{p^{1-\alpha}}{\sqrt{n}})+O(\frac{p^{1-\alpha}}{p^{\alpha/2}})$.

%$\sqrt(O(p^{1-\alpha})O(\frac{1}{n}+1/p^{\alpha/2}))=p^{\alpha/2-\frac{1}{2}}/{\sqrt{(n}}+p^{\alpha/4-\frac{1}{2}}=\frac{p^{\alpha/2}}{\sqrt{n}}+\frac{p^{\alpha/4}}{\sqrt{p}}$
%$=O({p^{\alpha/2}})O({\sqrt{\frac{\log(p)}{n}}})}$.

Applying the two proved theses of Theorem \ref{facPOET}
and Lemma \ref{Lemma11}, we can consequently prove, for each $k=1,\ldots,n$,
%Finally, we can show that
%$\max_{j\leq p, i\leq r}||\widehat{b}_j' \widehat{f}_i-b_j'f_i||=O(\frac{n^{1/4}}{p^{\alpha/2}}+\log(n)^{1/b_2}p^{1-\alpha}\sqrt{\frac{\log{p}}{n}})$, if $\alpha\leq 2/3+\frac{4}{3}\delta$, and
%$\max_{j\leq p, i\leq r}||\widehat{b}_j' \widehat{f}_i-b_j'f_i||=O(\frac{n^{1/4}}{p^{\alpha/2}}+\log(n)^{1/b_2}p^{\alpha/2-2\delta}\sqrt{\frac{\log{p}}{n}})$, if $2/3+\frac{4}{3}\delta<\alpha \leq 1$.
the third thesis:
\begin{equation}\max_{j\leq p, i\leq r}||\widehat{b}_j' \widehat{f}_k-b_j'f_k||=O\left(\frac{n^{1/4}p^{1-\alpha}}{p^{\alpha/2}}+\log(n)^{1/b_2}p^{\delta/2}\sqrt{\frac{\log{p}}{n}}\right)\end{equation}
simply following \cite{fan2013large}.

\subsection{Proof of Theorem \ref{facALCE}}\label{facMe}

Corollary 1 in \cite{farne2020large} prescribe that $\widehat{L}_{ALCE}=\widehat{U}_{ALCE}\widehat{D}_{ALCE}\widehat{U}_{ALCE}'$ is asymptotically consistent
if and only if $p^{\alpha+\delta}/n$ tends to $0$ as both $p$ and $n$ diverge.
Therefore, as these conditions are respected,
$\widehat{L}_{ALCE}$ and $\widehat{L}_{POET}$ both converge to $L^{*}$,
and the proof of Theorem \ref{facPOET} can be straightforwardly applied to UNALCE estimates
taking into account that now $\left(\frac{\widehat{D}_{ALCE}}{p}\right)^{-1}$
scale to $O(p^{1-\alpha-\delta})$ instead of $O(p^{1-\alpha})$.

\subsection{Proof of Theorem \ref{mine}}\label{proofs_3}

Let us decompose the minimization problems considered in Theorem \ref{mine}
conditioning on $Y_{pre}$ and $Z_{pre}$.
Suppose that the assumptions of Theorem \ref{mine} hold.
The problem in $L$ can be rewritten as
$$\min_{L \in \widehat{\mathcal{B}}(\widehat{r})} ||L-L^*||^2=\min_{L \in \widehat{\mathcal{B}}(\widehat{r})} ||L-Y_{pre}+Y_{pre}-L^*||^2 \leq
\min_{L \in \widehat{\mathcal{B}}(\widehat{r})} ||L-Y_{pre}||^2 + ||Y_{pre}-L^*||^2,$$
which is minimum if $L=\widehat{L}_{UNALCE}$ because of the optimal approximation property of principal components,
since $\widehat{L}_{UNALCE}$ is derived by the top $\widehat{r}$ principal components of ${Y}_{pre}$.

The problem in $S$ can be rewritten as follows. Suppose that, exploiting the assumption $diag(L)+diag(S)=diag(\widehat{\Sigma}_{ALCE})$,
we constrain ourselves within the class of matrices with diagonal $diag(\widehat{\Sigma}_{ALCE})-diag(L)$,
where $L$ has rank at most $r$ and belongs to $\widehat{\mathcal{B}}(\widehat{r})$. We call this space $\widehat{\mathcal{A}}_{diag}$.
%One way to see the problem is
%$$\min_{S \in \widehat{\Omega}} ||S-S^*||^2=\min_{S \in \widehat{\Omega}} ||S-S_{pre}+S_{pre}-S^*||^2 = ||\widehat{S}_{UNALCE}-S^*||^2.$$
%%because $S$ is uniquely determined by the assumption $diag(L)+diag(S)=diag(\widehat{\Sigma}_{ALCE})$.
%In alternative,
Defining ${\Sigma}_{pre}={Y}_{pre}+{Z}_{pre}$,
conditioning on $Y_{pre}$ and assuming the invariance of the off-diagonal elements in $\widehat{S}$,
we can write
$$\min_{S \in \widehat{\mathcal{A}}_{diag}} ||S-S^*||^2=\min_{L \in \widehat{\mathcal{B}}(\widehat{r})}||(\widehat{\Sigma}_{ALCE}-L)-(\Sigma_{pre}-Y_{pre})+(\Sigma_{pre}-Y_{pre})-(\Sigma^{*}-L^*)|| \leq$$%+(\Sigma_{pre}-Y_{pre})
$$\min_{L \in \widehat{\mathcal{B}}(\widehat{r})}||(\widehat{\Sigma}_{ALCE}-L)-(\Sigma_{pre}-Y_{pre})||+||(\Sigma_{pre}-Y_{pre})-(\Sigma^{*}-L^*)|| \leq $$
%$$\min_{L \in \widehat{\mathcal{B}}(\widehat{r})}||((\widehat{\Sigma}_{ALCE}-L)-(\Sigma_{pre}-Y_{pre})||+||(\Sigma_{pre}-Y_{pre})-(\Sigma^{*}-L^*)|| \leq$$
$$\min_{L \in \widehat{\mathcal{B}}(\widehat{r})}||(\widehat{\Sigma}_{ALCE}-\Sigma_{pre})||+||L-Y_{pre}||+||(\Sigma_{pre}-Y_{pre})-(\Sigma^{*}-L^*)||.$$
Since $(\widehat{\Sigma}_{ALCE}-\Sigma_{pre})$ and $(\Sigma_{pre}-Y_{pre})-(\Sigma^{*}-L^*)$ are fixed,
the entire problem boils down to $\min_{L \in \widehat{\mathcal{B}}(\widehat{r})} ||L-Y_{pre}||^2$, which is solved for $L=\widehat{L}_{UNALCE}$.
Therefore, we can write $\widehat{S}_{UNALCE}=\min_{S \in \widehat{\mathcal{A}}_{diag}} ||S-S^*||^2$, and the minimum amounts to $||\widehat{S}_{UNALCE}-S^*||^2$.

The problem in $\Sigma$ can be rewritten as
$$\min_{\Sigma \in \widehat{Y}} ||\Sigma-\Sigma^*||=\min_{\Sigma \in \widehat{Y}} ||L-L^*+S-S^*|| \leq$$
$$\leq \min_{L \in \widehat{\mathcal{B}}(\widehat{r})} ||L-L^*|| + \min_{S \in \widehat{\mathcal{A}}_{diag}} ||S-S^*|| \leq \min_{L \in \widehat{\mathcal{B}}(\widehat{r})} ||L-Y_{pre}|| + ||Y_{pre}-L^*||+||\widehat{S}_{UNALCE}-S^*||$$
or alternatively
$$\leq \min_{L \in \widehat{\mathcal{B}}(\widehat{r})} 2||L-Y_{pre}|| + ||Y_{pre}-L^*||+||(\widehat{\Sigma}_{ALCE}-\Sigma_{pre})||+||(\Sigma_{pre}-Y_{pre})-(\Sigma^{*}-L^*)||$$
which is minimum if $L=\widehat{L}_{UNALCE}$ and $S=\widehat{S}_{UNALCE}$.

The same optimality properties are transmitted to $S^{*-1}$ and $\Sigma^{*-1}$.
It is enough to recall that
$$||\widehat{S}^{-1}-S^{*-1}||\leq||{S}^{*-1}|| \times ||\widehat{S}-S^{*}|| \times ||{S}^{*-1}|| \leq 2\lambda_{min}(S^{*})||\widehat{S}-S^{*}||$$
and
$$||\widehat{\Sigma}^{-1}-\Sigma^{*-1}||\leq ||{\Sigma}^{*-1}|| \times ||\widehat{\Sigma}-\Sigma^{*}|| \times ||{\Sigma}^{*-1}||\leq 2\lambda_{min}(\Sigma^{*})||\widehat{\Sigma}-\Sigma^{*}||,$$
such that $\widehat{S}_{UNALCE}^{-1}$ and $\widehat{\Sigma}_{UNALCE}^{-1}$ minimize $||\widehat{S}^{-1}-S^{*-1}||$ and $||\widehat{\Sigma}^{-1}-\Sigma^{*-1}||$ respectively under the assumptions of Theorem \ref{mine}.

\subsection{Proof of Corollary \ref{Coroll_pl}}
In order to prove the first part of the corollary, we can state that
$$||\mathbb{P}_{T'^\perp} (\widehat{L}_{UNALCE}-L^*)||\leq ||(\widehat{L}_{UNALCE}-L^*)||\leq$$
$$||(\widehat{L}_{UNALCE}-\widehat{L}_{ALCE}+\widehat{L}_{ALCE}-L^*)||\leq$$ $$||(\widehat{L}_{ALCE}-L^*)||+||(\widehat{L}_{UNALCE}-\widehat{L}_{ALCE})|| \leq (C+1)\psi,$$
where $C\geq \xi(T)$,
because $\widehat{L}_{UNALCE}=\widehat{U}_{ALCE} (\widehat{D}_{ALCE} +\lambda I_r)\widehat{U}_{ALCE} '=\widehat{L}_{ALCE}+\widehat{U}_{ALCE} \lambda I_r \widehat{U}_{ALCE} '$. %and
Relying on $\widehat{S}_{UNALCE}=\widehat{S}_{ALCE}-diag(\widehat{U}_{ALCE} (\widehat{D}_{ALCE} +\lambda I_r)\widehat{U}_{ALCE} )$, we then have
$$g_\gamma(\widehat{\Sigma}_{UNALCE}-\Sigma^{*},\widehat{\Sigma}_{UNALCE}-\Sigma^{*})=
||\widehat{\Sigma}_{UNALCE}-\widehat{\Sigma}_{ALCE}
+\widehat{\Sigma}_{ALCE}-\Sigma^{*}||\leq$$
$$||\widehat{\Sigma}_{UNALCE}-\widehat{\Sigma}_{ALCE}||+||\widehat{\Sigma}_{ALCE}-\Sigma^{*}||\leq$$
$$||\widehat{L}_{UNALCE}-\widehat{L}_{ALCE}||+||\widehat{S}_{UNALCE}-\widehat{S}_{ALCE}||+||\widehat{\Sigma}_{ALCE}-\Sigma^{*}||=$$
$$|\widehat{U}_{ALCE} (\lambda I_r)\widehat{U}_{ALCE}||+||-\widehat{U}_{ALCE} (\lambda I_r)\widehat{U}_{ALCE}||+||\widehat{\Sigma}_{ALCE}-\Sigma^{*}||=$$
$$(C+2)\lambda,$$
with $C\leq 11$ (see \cite{luo2011high}).

%Conditionally on $\widehat{\Sigma}_{ALCE}=\widehat{L}_{ALCE,PCA}+\widehat{S}_{ALCE,PCA}$, we have
In order to prove the second part of the corollary, we recall from \cite{farne2020large} that $||\widehat{L}_{ALCE}-Y_{pre}||-||\widehat{L}_{UNALCE}-Y_{pre}||\leq \psi$ and
$||\widehat{\Sigma}_{ALCE}-\Sigma_{pre}||-||\widehat{\Sigma}_{UNALCE}-\Sigma_{pre}|| \leq \psi$.
%because it holds $||\mathbb{P}_{T'^\perp} (L-L^*)||\leq ||\mathbb{P}_{T'^\perp} (L-\widehat{L}_{ALCE,PCA})||+||\mathbb{P}_{T'^\perp} (\widehat{L}_{ALCE}-L^*)||$,
%which is minimum for $L=\widehat{L}_{UNALCE}$ because $||\mathbb{P}_{T'^\perp} (L-\widehat{L}_{ALCE,PCA})||=0$.
Then we can write

$$||\mathbb{P}_{T'^\perp} (\widehat{L}_{ALCE}-Y_{pre}+Y_{pre}-L^*)||\leq$$
$$||(\widehat{L}_{ALCE}-Y_{pre}+Y_{pre}-L^*)||\leq$$ $$||\widehat{L}_{UNALCE}-\widehat{L}_{ALCE}||+||(\widehat{L}_{UNALCE}-Y_{pre})||+||(Y_{pre}-L^{*})||$$
Therefore, $$||\mathbb{P}_{T'^\perp} (\widehat{L}_{ALCE}-L^*)|-||\mathbb{P}_{T'^\perp} (\widehat{L}_{UNALCE}-L^*)||
\leq$$
$$||(\widehat{L}_{ALCE}-L^*)|| - ||(\widehat{L}_{UNALCE}-L^*)|| \leq\psi.$$
%$C+1\leq \xi(T)$.
%(see \cite{farne2020large}).

Similarly, conditioning on $\Sigma_{pre}$ and recalling \cite{farne2020large}, we have
$$g_\gamma(\widehat{\Sigma}_{ALCE}-\Sigma^{*},\widehat{\Sigma}_{ALCE}-\Sigma^{*})-
g_\gamma(\widehat{\Sigma}_{UNALCE}-\Sigma^{*},\widehat{\Sigma}_{UNALCE}-\Sigma^{*})\leq$$
$$||\widehat{\Sigma}_{ALCE}-\Sigma^{*}||%-\widehat{\Sigma}_{UNALCE}
-||\widehat{\Sigma}_{UNALCE}-\Sigma^{*}||\leq \psi.$$

\subsection{Proof of Theorem \ref{eigen}}\label{proofs_4}

Restricting to $\widehat{\mathcal{Y}}=\widehat{\mathcal{B}}(\widehat{r})+\widehat{\mathcal{A}}(\widehat{s})$ and conditioning on $Y_{pre}$ and $Z_{pre}$
(which in turn rely on $\Sigma_n$), we have proved in Theorem \ref{mine} that
$\widehat{L}_{UNALCE}=\min_{L \in \widehat{\mathcal{B}}(\widehat{r})} ||L-L^*||^2$, $\widehat{S}_{UNALCE}=\min_{S \in \widehat{\mathcal{A}}(\widehat{s})} ||S-S^*||^2$
and $\widehat{\Sigma}_{UNALCE}=\min_{\Sigma \in \widehat{Y}} ||\Sigma-\Sigma^*||^2$.
According to \cite{ledoit2004well}, we can write
\begin{eqnarray}
\frac{1}{p} E\left[\sum_{i=1}^p (\widehat{\lambda}_{L,i}-\mu_{L})^2 \vert \Sigma_n \right]=\frac{1}{p}\sum_{i=1}^p(\lambda_{L,i}-\mu_{L})^2+E(||\widehat{L}-L^{*}||^2 \vert \Sigma_n),
\label{eigL}\\
\frac{1}{p} E\left[\sum_{i=1}^p (\widehat{\lambda}_{S,i}-\mu_{S})^2 \vert \Sigma_n \right]=\frac{1}{p}\sum_{i=1}^p(\lambda_{S,i}-\mu_{S})^2+E(||\widehat{S}-S^{*}||^2 \vert \Sigma_n),
\label{eigS}
\\
\frac{1}{p} E\left[\sum_{i=1}^p (\widehat{\lambda}_{\Sigma,i}-\mu_{\Sigma})^2 \vert \Sigma_n \right]=\frac{1}{p}\sum_{i=1}^p(\lambda_{\Sigma,i}-\mu_{L})^2+E(||\widehat{\Sigma}-\Sigma^{*}||^2 \vert \Sigma_n),
\label{eigSigma}
\end{eqnarray}
where $\mu_{L}=tr({L}^{*})/p$, $\mu_{S}=tr({S}^{*})/p$  and $\mu_{\Sigma}=tr({\Sigma}^{*})/p$.
%Suppose that $\psi$ is sufficiently small, %\rightarrow 0$, %i.e. $\frac{{p}^{\alpha}}{\sqrt{n}}$ tends to $0$ while $n=o(p^{2\alpha})$.
%such that $E(trace(\widehat{\lambda}_{L})/p)\sim tr(L^{*})/p$ and $E(trace(\widehat{\lambda}_{\Sigma})/p)\sim tr(\Sigma^{*})/p$.
As a consequence, since $E(||\widehat{L}~-~L^{*}||^2 \vert \Sigma_n)$, $E(||\widehat{S}-S^{*}||^2 \vert \Sigma_n)$ and
$E(||\widehat{\Sigma}-\Sigma^{*}||^2 \vert \Sigma_n)$ are minimum under the assumptions of Theorem \ref{mine},
we can write $\widehat{L}_{UNALCE}=\min_{L \in \widehat{\mathcal{B}}(\widehat{r})}\frac{1}{p} E\left[\sum_{i=1}^p (\widehat{\lambda}_{L,i}-\mu_{L})^2\right]$,
$\widehat{S}_{UNALCE}=\min_{S \in \widehat{\mathcal{A}}_{diag}}\frac{1}{p} E\left[\sum_{i=1}^p (\widehat{\lambda}_{S,i}-\mu_{S})^2\right]$
and $\widehat{\Sigma}_{UNALCE}=\min_{\Sigma \in \widehat{Y}}\frac{1}{p} E\left[\sum_{i=1}^p (\widehat{\lambda}_{\Sigma,i}-\mu_{\Sigma})^2\right]$.
These results mean that within the classes of algebraically consistent estimates
the eigenvalues of $\widehat{L}_{UNALCE}$, $\widehat{S}_{UNALCE}$, and $\widehat{\Sigma}_{UNALCE}$
are the most concentrated possible around their respective true means.

The optimality properties of the eigenvalues of $S^{*}$ and $\Sigma^{*}$ estimated by UNALCE
are transmitted to $S^{*-1}$ and $\Sigma^{*-1}$. In fact, we know that $E(||\widehat{S}^{-1}~-~S^{*-1}||^2 \vert \Sigma_n)$ and $E(||\widehat{\Sigma}^{-1}~-~\Sigma^{*-1}||^2 \vert \Sigma_n)$ are the minimum possible under the assumptions of Theorem \ref{mine}.
%because in \cite{farne2017thesis} we need $\phi^{-1}\geq 2\lambda_{min}^{-1}$, i.e., $\phi \leq \frac{1}{2} \lambda_{min}$, in order to achieve $||\widehat{\Sigma}^{-1}-\Sigma^{-1}||\leq\phi$, such that the bound above is consequently necessary.
Since it also holds
$$
\frac{1}{p} E\left[\sum_{i=1}^p (\widehat{\lambda}_{S^{-1},i}-\mu_{S^{-1}})^2 \vert \Sigma_n \right]=\frac{1}{p}\sum_{i=1}^p(\lambda_{S^{-1},i}-\mu_{S^{-1}})^2+E(||\widehat{S}^{-1}-S^{*-1}||^2 \vert \Sigma_n),
$$
$$
\frac{1}{p} E\left[\sum_{i=1}^p (\widehat{\lambda}_{\Sigma^{-1},i}-\mu_{\Sigma^{-1}})^2 \vert \Sigma_n \right]=\frac{1}{p}\sum_{i=1}^p(\lambda_{\Sigma^{-1},i}-\mu_{\Sigma^{-1}})^2+E(||\widehat{\Sigma}^{-1}-\Sigma^{*-1}||^2 \vert \Sigma_n),
$$
%if $||\widehat{\Sigma}-\Sigma^{*}||$ is minimum, $||\widehat{\Sigma}^{*-1}-\Sigma^{*-1}||$ is also in turn,
%and consequently the equation
%$$
%\frac{1}{p} E\left[\sum_{i=1}^p (\widehat{\lambda}_{\Sigma^{-1},i}-\mu_{\Sigma^{-1}})^2 \vert \Sigma_n \right]=\frac{1}{p}\sum_{i=1}^p(\lambda_{\Sigma^{-1},i}-\mu_{\Lambda^{-1}})^2+E(||\widehat{\Sigma}^{-1}-\Sigma^{*-1}||^2 \vert \Sigma_n),
%$$
where $\mu_{S^{-1}}$ and $\mu_{\Sigma^{-1}}$ are the mean eigenvalues of $S^{*-1}$ and $\Sigma^{*-1}$,
we are allowed to conclude that $\widehat{S}^{-1}_{UNALCE}=\min_{S \in \widehat{\mathcal{A}}_{diag}} \sum_{i=1}^p \widehat{S}^{-1}_{UNALCE}=(\widehat{\lambda}_{S^{-1},i}-\mu_{S^{-1}})^2$ and $\widehat{\Sigma}^{-1}_{UNALCE}=\min_{\Sigma \in \widehat{Y}} \sum_{i=1}^p (\widehat{\lambda}_{\Sigma^{-1},i}-\mu_{\Sigma^{-1}})^2$ .
%This allows to use the result about the optimality of the spectrum of $\widehat{\Sigma^{-1}}_{UNALCE}$ and the asymptotic results about the moments of the ESD.

\subsection{Proof of Corollary \ref{esd}}\label{proof_esd}
We show the proof with respect to $\widehat{L}_{UNALCE}$ as the extension to $\widehat{S}_{UNALCE}$, $\widehat{\Sigma}_{UNALCE}$,
$\widehat{S}^{-1}_{UNALCE}$, $\widehat{\Sigma}^{-1}_{UNALCE}$ is straightforward.
From Theorem \ref{eigen} we know that
\begin{eqnarray}
sum_L=\sum_{i=1}^p (\widehat{\lambda}_{\widehat{L}_{UNALCE},i}-\mu_{L})^2= \nonumber\\
=\sum_{i=1}^p \widehat{\lambda}_{\widehat{L}_{UNALCE},i}^2 + p \mu_{L}^2 -2 \mu_L \sum_{i=1}^p \widehat{\lambda}_{\widehat{L}_{UNALCE},i} \nonumber
\end{eqnarray}
is minimum into the recovered low rank matrix variety.
Then we note that $sum_L/p$ can be rewritten as
$\frac{1}{p} tr(\widehat{L}_{UNALCE}-\mu_L I_p)^2$.
We know that $\frac{1}{p} (tr(\widehat{L}_{UNALCE}-\mu_L I_p)^2)$ is the second moment of $\rho(z)_{\widehat{L}_{UNALCE}-{L^*}}$,
because $tr(\mu_L I_p)=tr(L^{*})$ and
$tr(\widehat{L}_{UNALCE}-\mu_L I_p)^2=tr(\widehat{L}_{UNALCE}-L^{*})^2$.
Therefore, the claim on the second moment of $\rho(z)_{\widehat{L}_{UNALCE}-L^{*}}^p$ is proved.
Concerning the claim on the first moment of $\rho(z)_{\widehat{L}_{UNALCE}-L^{*}}$,
it is sufficient to note that $tr(\widehat{L}_{UNALCE}-\mu_L I_p)=tr(\widehat{L}_{UNALCE}-L^{*})$ tends to $0$
as $\psi=\frac{1}{\xi(T)} \frac{p^{\alpha}}{\sqrt{n}}$ tends to $0$.

\subsection{Proof of Corollary \ref{Coroll_B}}

Since the eigenvalues of $\widehat{L}_{UNALCE}$ are the most concentrated around their mean under the assumptions of Theorem \ref{mine}, the same holds for $\widehat{B}_{UNALCE}$, because its eigenvalues are the square root of the ones of $\widehat{L}_{UNALCE}$ and the variance is a monotonic operator. Therefore, according to \cite{ledoit2004well},
$\min_{B,{L}={B}{B}'\in\widehat{\mathcal{B}}} ||\widehat{B}-B|| \vert \Sigma_n$ is solved by $\widehat{B}=\widehat{B}_{UNALCE}$ under the constraint $\widehat{B}'\widehat{B}$ diagonal and $\sum_{i=1}^p ||\widehat{b}_i||=\max$.

\subsection{Proof of Theorem \ref{bartlett_opt}}\label{proof_bart}

We start considering the loss $||\widehat{f}_{k,B}-f_{k,B}||$.
Since $\widehat{f}_{k,B}=(\widehat{B}'\widehat{S}^{-1}\widehat{B})^{-1}\widehat{B}'\widehat{S}^{-1} x_k$, $k=1,\ldots,n$,
that loss is majorized by $||(\widehat{B}'\widehat{S}^{-1}\widehat{B})^{-1}\widehat{B}'\widehat{S}^{-1}-(B'S^{*-1}B)^{-1}B'S^{*-1}|| \times ||x_k||$.
%The $r$ eigenvalues of $(\widehat{B}'\widehat{S}^{-1}\widehat{B})^{-1}\widehat{B}'\widehat{S}^{-1}$ are the same as the ones of $(\widehat{S}^{-1}\widehat{B})(\widehat{B}'\widehat{S}^{-1}\widehat{B})^{-1}\widehat{B}'$.
The eigenvalues of $\widehat{B}'\widehat{S}^{-1}\widehat{B}$
coincide with the ones of $\widehat{S}^{-1}\widehat{B}\widehat{B}'=\widehat{S}^{-1}\widehat{L}$. According to \cite{ledoit2004well}, the expected variance of the estimated eigenvalues around their true mean depends on the true variance and the squared bias in spectral norm of the overall estimate. Therefore, conditioning on $\Sigma_n$, we can focus on the spectral losses $||\widehat{B}' \widehat{S}^{-1}-B' S^{*-1}||$ and $||\widehat{S}^{-1}\widehat{L}-S^{*-1}L^{*}||$. % because the spectral loss $||\widehat{B}'-B'||$
%is minimum for $\widehat{B}=\widehat{B}_{UNALCE}$ for Corollary \ref{Coroll_B}.

%is $||\widehat{B}\widehat{B}'\widehat{S}^{-1}x_k-BB'S^{*-1}x_k||$, which
%entirely depends on $||\widehat{B}\widehat{B}'\widehat{S}^{-1}-BB'S^{*-1}||=||\widehat{L}'\widehat{S}^{-1}-L^{*}S^{*-1}||$.
We first focus on $||\widehat{S}^{-1}\widehat{L}-S^{*-1}L^{*}||$.
Conditioning on $Y_{pre}$ and $Z_{pre}$, we can write
$$||\widehat{S}^{-1}\widehat{L}-S^{*-1}L^{*}|| \leq ||\widehat{S}^{-1}\widehat{L}-Z_{pre}^{-1}Y_{pre}|| +
||Z_{pre}^{-1}Y_{pre}-S^{*-1}L^{*}||,$$ where the term $||Z_{pre}^{-1}Y_{pre}-S^{*-1}L^{*}||$
entirely depends on $\Sigma_n$.

We now consider two generic estimates $L$ and $S$. Under the assumptions of Theorem \ref{mine}, we must constrain our search by setting $off-diag(S)=off-diag(\widehat{S}_{ALCE})$ and $diag(S)=diag(\widehat{\Sigma}_{ALCE}-L)$, $L \in \mathcal{L}(\widehat{r})$.
%We call this set $\mathcal{S}_{diag}(\widehat{s})$.
Therefore, conditioning on ${\Sigma}_{pre}={Y}_{pre}+{Z}_{pre}$, we can write
$$||\widehat{S}^{-1}\widehat{L}-Z_{pre}^{-1}Y_{pre}||=||({\Sigma}_{pre}-L)^{-1}L-({\Sigma}_{pre}-Y_{pre})^{-1}Y_{pre}||,$$
%$||\widehat{B}\widehat{B}'\widehat{S}^{-1}-BB'\Sigma^{*-1}||$
with $off-diag(L)=off-diag({Y}_{pre})$.
We apply the formula $({\Sigma}_{pre}-L)^{-1}=\sum_{k=0}^{\infty} ({\Sigma}_{pre}^{-1}L)^{k} {\Sigma}_{pre}^{-1}$,
which leads to $$||({\Sigma}_{pre}-L)^{-1}L-({\Sigma}_{pre}-Y_{pre})^{-1}Y_{pre}||= \left \vert \left \vert \sum_{k=0}^{\infty} ({\Sigma}_{pre}^{-1}L)^{k} {\Sigma}_{pre}^{-1}L-\sum_{k=0}^{\infty} ({\Sigma}_{pre}^{-1}Y_{pre})^{k} {\Sigma}_{pre}^{-1}Y_{pre}
\right \vert \right \vert.$$

Conditioning on $Y_{pre}$, we can write $L=Y_{pre}+\Delta_{L,pre}$.
Therefore, it follows
$$||({\Sigma}_{pre}-L)^{-1}L-({\Sigma}_{pre}-Y_{pre})^{-1}Y_{pre}||= \left \vert \left \vert
\sum_{k=0}^{\infty} \Sigma_{pre}^{-k} (L^{k}-Y_{pre}^{k}) {\Sigma}_{pre}^{-1} \Delta_{L,pre}
\right \vert \right \vert.
$$

Applying Cauchy-Schwartz inequality, we obtain
$$\left \vert \left \vert  ({\Sigma}_{pre}-L)^{-1}L-({\Sigma}_{pre}-Y_{pre})^{-1}Y_{pre}\right \vert \right \vert\leq$$
\begin{equation}\sum_{k=0}^\infty ||{\Sigma}_{pre}^{-k}||\times||L^{k}-Y_{pre}^{k}|| \times ||{\Sigma}_{pre}^{-1}|| \times ||\Delta_{L,pre}|| \label{subth1}\end{equation}

Recalling Theorem \ref{eigen}, which states that the variance of the eigenvalues of $\widehat{L}_{UNALCE}$
are the most concentrated possible around their true mean within the recovered low rank variety,
we note that this holds for any power of $\widehat{L}_{UNALCE}$,
$\widehat{L}_{UNALCE}^{k}$ (with $k \ne 0$), due to the monotonicity of the variance operator.
For this reason, $||L^{k}-Y_{pre}^{k}||$ is the minimum possible for $L=\widehat{L}_{UNALCE}$ at any $k$.
Under the assumptions of Theorem \ref{mine}, the same holds for $||\Delta_{L,pre}||$, and
the problem in $S$ is solved by $\widehat{S}_{UNALCE}$.
Therefore, $||\widehat{S}^{-1}\widehat{L}-Z_{pre}^{-1}Y_{pre}||$ is minimum for $\widehat{L}_{UNALCE}$ and $\widehat{S}_{UNALCE}$.

Conditioning on $Y_{pre}$ and $Z_{pre}$, we can write
$$||\widehat{B}'\widehat{S}^{-1}-B' S^{*-1}|| \leq ||\widehat{B}'\widehat{S}^{-1}-B_{pre}' Z_{pre}^{-1}|| +
||B_{pre}' Z_{pre}^{-1}-B' S^{*-1}||,$$ where the term $||B_{pre}' Z_{pre}^{-1}-B' S^{*-1}||$
entirely depends on $\Sigma_n$.
Conditioning on $B_{pre}$, obtained defining $Y_{pre}=B_{pre}B_{pre}'$,
applying the same framework we can write $B=B_{pre}+\Delta_{B,pre}$, and
\begin{equation} ||\widehat{B}'\widehat{S}^{-1}-B_{pre}' Z_{pre}^{-1}|| \leq ||\Delta_{B,pre}|| \times \sum_{k=0}^\infty ||{\Sigma}_{pre}^{-k}||\times||(L^{k}-Y_{pre}^{k})|| \times ||{\Sigma}_{pre}^{-1}||, \label{subth2} \end{equation}
which is minimum for $\widehat{B}=\widehat{B}_{UNALCE}$ from Theorem \ref{mine}.

Starting from
$$||(\widehat{B}'\widehat{S}^{-1}\widehat{B})^{-1}\widehat{B}'\widehat{S}^{-1}-(B'S^{*-1}B)^{-1}B'S^{*-1}||\leq$$
$$||(\widehat{B}'\widehat{S}^{-1}\widehat{B})^{-1}\widehat{B}'\widehat{S}^{-1}-(B_{pre}'Z_{pre}^{-1}B_{pre})^{-1}B_{pre}'Z_{pre}^{-1}||+$$
$$+||(B_{pre}'Z_{pre}^{-1}B_{pre})^{-1}B_{pre}'Z_{pre}^{-1}-(B'S^{*-1}B)^{-1}B'S^{*-1}||,$$
noting that
$$(\widehat{B}'\widehat{S}^{-1}\widehat{B})^{-1}\widehat{B}'\widehat{S}^{-1}-(B_{pre}'Z_{pre}^{-1}B_{pre})^{-1}B_{pre}'Z_{pre}^{-1}=$$
$$=[(\widehat{B}'\widehat{S}^{-1}\widehat{B})^{-1}-(B_{pre}'Z_{pre}^{-1}B_{pre})^{-1}][(\widehat{B}'\widehat{S}^{-1}+B_{pre}'Z_{pre}^{-1})]+$$
$$-(\widehat{B}'\widehat{S}^{-1}\widehat{B})^{-1}B_{pre}'Z_{pre}^{-1}+(B_{pre}'Z_{pre}^{-1}B_{pre})^{-1}\widehat{B}'\widehat{S}^{-1},$$
we obtain
$$||(\widehat{B}'\widehat{S}^{-1}\widehat{B})^{-1}\widehat{B}'\widehat{S}^{-1}-(B_{pre}'Z_{pre}^{-1}B_{pre})^{-1}B_{pre}'Z_{pre}^{-1}||\leq$$
$$\leq ||[(\widehat{B}'\widehat{S}^{-1}\widehat{B})^{-1}-(B_{pre}'Z_{pre}^{-1}B_{pre})^{-1}]||\times||[(\widehat{B}'\widehat{S}^{-1}+B_{pre}'Z_{pre}^{-1})]||\times$$
$$\times||(\widehat{B}'\widehat{S}^{-1}\widehat{B})^{-1}B_{pre}'Z_{pre}^{-1}|| \times ||(B_{pre}'Z_{pre}^{-1}B_{pre})^{-1}\widehat{B}'\widehat{S}^{-1}||.\label{rel}$$

Conditioning on $B_{pre}$ and $Z_{pre}$,
$$||[(\widehat{B}'\widehat{S}^{-1}+B_{pre}'Z_{pre}^{-1})]||=||[(\widehat{B}'\widehat{S}^{-1}-B_{pre}'Z_{pre}^{-1}+2B_{pre}'Z_{pre}^{-1}]||\leq$$
$$ \leq ||\widehat{B}'\widehat{S}^{-1}-B_{pre}'Z_{pre}^{-1}||+2||B_{pre}'Z_{pre}^{-1}||$$
Therefore, the first and the second multiplicative factors are minimum for $\widehat{B}=\widehat{B}_{UNALCE}$ and $\widehat{S}=\widehat{S}_{UNALCE}$
for \ref{subth1} and \ref{subth2} respectively.

Concerning $||(\widehat{B}'\widehat{S}^{-1}\widehat{B})^{-1}B_{pre}'Z_{pre}^{-1}||$ and $||(B_{pre}'Z_{pre}^{-1}B_{pre})^{-1}\widehat{B}'\widehat{S}^{-1}||$, it is sufficient to write
$(\widehat{B}'\widehat{S}^{-1}\widehat{B})^{-1}=(B_{pre}'Z_{pre}^{-1}B_{pre})^{-1}+\Delta_{B'S^{-1}B}$
and
$\widehat{B}'\widehat{S}^{-1}=B_{pre}'Z_{pre}^{-1}+\Delta_{B'S^{-1}}$,
%and for \ref{subth1} and \ref{subth2}
such that both norms are minimum for $B=\widehat{B}_{UNALCE}$ and $S=\widehat{S}_{UNALCE}$
because $||\Delta_{B'S^{-1}B}||$ and $||\Delta_{B'S^{-1}}||$ are minimum for \ref{subth1}
and \ref{subth2}.

Finally, we can extend the validity to $\widehat{B}\widehat{f}_{k,B}-Bf_{k,B}$ noting that
$||\widehat{B}\widehat{f}_{k,B}-Bf_{k,B}||\leq ||\widehat{B}(\widehat{B}'\widehat{S}^{-1}\widehat{B})^{-1}\widehat{B}'\widehat{S}^{-1}-B(B'S^{*-1}B)^{-1}B'S^{*-1}|| \times ||x_k||$ and claiming that
\begin{eqnarray}
%$(\widehat{B}'(\widehat{S^})^{-1}\widehat{B})^{-1}\widehat{B}'\widehat{S}^{-1}\widehat{B}$
%\widehat{B}\widehat{f}_{k,B}-Bf_{k,B}=
\widehat{B}(\widehat{B}'\widehat{S}^{-1}\widehat{B})^{-1}\widehat{B}'\widehat{S}^{-1}-B(B'S^{*-1}B)^{-1}B'S^{*-1}= \nonumber \\
=\widehat{B}(\widehat{B}'\widehat{S}^{-1}\widehat{B})^{-1}\widehat{B}'\widehat{S}^{-1}-B(B'S^{*-1}B)^{-1}B'S^{*-1}+ \nonumber \\ %$ (||A||=||A'||)
+B(\widehat{B}'\widehat{S}^{-1}\widehat{B})^{-1}\widehat{B}'\widehat{S}^{-1}-B(\widehat{B}'\widehat{S}^{-1}\widehat{B})^{-1}\widehat{B}'\widehat{S}^{-1}+ \nonumber \\
+(\widehat{B}-B)(B'S^{*-1}B)^{-1}B'S^{*-1}-(\widehat{B}-B)(B'S^{*-1}B)^{-1}B'S^{*-1} \nonumber
\end{eqnarray}
which becomes
\begin{eqnarray}
(\widehat{B}-B)(\widehat{B}'\widehat{S}^{-1}\widehat{B})^{-1}\widehat{B}'\widehat{S}^{-1}+B((\widehat{B}'\widehat{S}^{-1}\widehat{B})^{-1}\widehat{B}'
\widehat{S}^{-1}-(B'S^{*-1}B)^{-1}B'S^{*-1})\nonumber \\
-(\widehat{B}-B)(B'S^{*-1}B)^{-1}B'S^{*-1}+(\widehat{B}-B)(B'S^{*-1}B)^{-1}B'S^{*-1}= \nonumber \\
(\widehat{B}-B) ((\widehat{B}'\widehat{S}^{-1}\widehat{B})^{-1}\widehat{B}'\widehat{S}^{-1}-(B'S^{*-1}B)^{-1}B'S^{*-1}))\nonumber \\
+B((\widehat{B}'\widehat{S}^{-1}\widehat{B})^{-1}\widehat{B}'\widehat{S}^{-1}-(B'S^{*-1}B)^{-1}B'S^{*-1}))+(\widehat{B}-B)(B'S^{*-1}B)^{-1}B'S^{*-1}. \nonumber
\end{eqnarray}
Since $(\widehat{B}-B)$ and $(\widehat{B}'\widehat{S}^{-1}\widehat{B})^{-1}\widehat{B}'\widehat{S}^{-1}-(B'S^{*-1}B)^{-1}B'S^{*-1}$
are minimum for $\widehat{B}=\widehat{B}_{UNALCE}$ and $\widehat{S}=\widehat{S}_{UNALCE}$ as we previously proved,
we can derive that also $\widehat{B}(\widehat{B}'\widehat{S}^{-1}\widehat{B})^{-1}\widehat{B}'\widehat{S}^{-1}-B(B'S^{*-1}B)^{-1}B'S^{*-1}=\widehat{B}\widehat{f}_{B}-Bf$
is the minimum possible for the same matrices, thus proving the thesis.

\subsection{Proof of Theorem \ref{thompson_opt}}\label{proof_thom}

We start considering the loss $||\widehat{B}\widehat{f}_{k,T}-Bf_{k,T}||$. %and $||\widehat{B}\widehat{f}_{k,T}~-~Bf||$.
For the definition of $\widehat{f}_{k,T}$, that loss is majorized by
$||\widehat{B}\widehat{B}'\widehat{\Sigma}^{-1}-BB'\Sigma^{*-1}||\times ||x_k||=||\widehat{L}\widehat{\Sigma}^{-1}-L\Sigma^{*-1}||\times||x_k||$.
%$||\widehat{B}_{UN,T}\widehat{B}_{UN,T}'\widehat{\Sigma}_{UN}^{-1}-BB'\Sigma^{*-1}||$.\\
Conditioning on $Y_{pre}$, $Z_{pre}$ and $\Sigma_{pre}$, we write
\begin{equation}||\widehat{L}\widehat{\Sigma}^{-1}-L^{*}\Sigma^{*-1}|| \leq ||\widehat{L}\widehat{\Sigma}^{-1}-Y_{pre}\Sigma_{pre}^{-1}|| + ||Y_{pre}\Sigma_{pre}^{-1}-L^{*}\Sigma^{*-1}||.\label{inverse_sum}\end{equation}

Then, we apply the following formula for the inverse of a sum: $$\widehat{\Sigma}^{-1}=(\widehat{S}+\widehat{L})^{-1}=\widehat{S}^{-1}-\widehat{S}^{-1}(I_p+\widehat{L}\widehat{S}^{-1})^{-1}\widehat{L}\widehat{S}^{-1}.$$
We observe that $||I_p+\widehat{L}\widehat{S}^{-1}-(I_p-Y_{pre} Z_{pre}^{-1})||=||\widehat{L}\widehat{S}^{-1}-Y_{pre} Z_{pre}^{-1}||$,
such that $(I_p+\widehat{L}_{UNALCE}\widehat{S}_{UNALCE}^{-1})$ inherits the optimality properties of $\widehat{L}_{UNALCE}\widehat{S}_{UNALCE}^{-1}$ previously proved. In addition, the same optimality property is transmitted to $(I_p+\widehat{L}_{UNALCE}\widehat{S}_{UNALCE}^{-1})^{-1}$,
for the consequences of Theorems \ref{mine} and \ref{eigen}.

Therefore,
$$||\widehat{\Sigma}^{-1}-\Sigma_{pre}^{-1}||=$$
$$||\widehat{S}^{-1}-\widehat{S}^{-1}(I_p+\widehat{L}\widehat{S}^{-1})^{-1}\widehat{L}\widehat{S}^{-1}-[Z_{pre}^{-1}-Z_{pre}^{-1}(I_p+Y_{pre}Z_{pre}^{-1})^{-1}Y_{pre}Z_{pre}^{-1}]||\leq$$
$$||\widehat{S}^{-1}-Z_{pre}^{-1}||+||\widehat{S}^{-1}(I_p+\widehat{L}\widehat{S}^{-1})^{-1}\widehat{L}\widehat{S}^{-1}-Z_{pre}^{-1}(I_p+Y_{pre}Z_{pre}^{-1})^{-1}Y_{pre}Z_{pre}^{-1}||.$$

The first term is minimum for $\widehat{S}_{UNALCE}$ from Theorem \ref{mine}.
%Conditioning on $\Sigma_{pre}$, $\Sigma_{pre}=Y_{pre}+\Delta_{L,pre}$,
%the second term can be factorized as follows.
Noting that
$$\widehat{S}^{-1}(I_p+\widehat{L}\widehat{S}^{-1})^{-1}\widehat{L}\widehat{S}^{-1}-Z_{pre}^{-1}(I_p+Y_{pre}Z_{pre}^{-1})^{-1}Y_{pre}Z_{pre}^{-1}=$$
$$=[\widehat{S}^{-1}(I_p+\widehat{L}\widehat{S}^{-1})^{-1}-Z_{pre}^{-1}(I_p+Y_{pre}Z_{pre}^{-1})^{-1}][\widehat{L}\widehat{S}^{-1}+Y_{pre}Z_{pre}^{-1}]+$$
$$+Z_{pre}^{-1}(I_p+Y_{pre}Z_{pre}^{-1})\widehat{L}\widehat{S}^{-1}-\widehat{S}^{-1}(I_p+\widehat{L}\widehat{S}^{-1})^{-1}Y_{pre}Z_{pre}^{-1},$$
and in turn
$$[\widehat{S}^{-1}(I_p+\widehat{L}\widehat{S}^{-1})^{-1}-Z_{pre}^{-1}(I_p+Y_{pre}Z_{pre}^{-1})^{-1}]=$$
$$(\widehat{S}^{-1}-{Z}_{pre}^{-1})[(I_p+\widehat{L}\widehat{S}^{-1})^{-1}+(I_p+Y_{pre}Z_{pre}^{-1})^{-1}]+$$
$$-\widehat{S}^{-1}(I_p+Y_{pre}Z_{pre}^{-1})^{-1}+Z_{pre}^{-1}(I_p+\widehat{L}\widehat{S}^{-1})^{-1},$$
conditioning on $Y_{pre}$, $Z_{pre}$ and $\Sigma_{pre}$,
it follows from Theorem \ref{mine}
that the minimum is for $\widehat{L}=\widehat{L}_{UNALCE}$ and $\widehat{S}=\widehat{S}_{UNALCE}$.
%$||\widehat{L}\widehat{\Sigma}^{-1}-Y_{pre}\Sigma_{pre}^{-1}||$

Then, moving from \ref{inverse_sum} and applying Cauchy-Scwartz inequality, we obtain
$$||\widehat{L}\widehat{\Sigma}^{-1}-Y_{pre}\Sigma_{pre}^{-1}||\leq
||\widehat{L}-Y_{pre}||\times||\widehat{\Sigma}^{-1}+\Sigma_{pre}^{-1}||
\times||\widehat{L}\Sigma_{pre}^{-1}|| \times ||Y_{pre}\widehat{\Sigma}^{-1}||,$$
and conditioning on $Y_{pre}$ and $\Sigma_{pre}$,
since $$||\widehat{\Sigma}^{-1}+\Sigma_{pre}^{-1}|| \leq ||\widehat{\Sigma}^{-1}-\Sigma_{pre}^{-1}||+2||\Sigma_{pre}^{-1}||,$$
$\widehat{L}=Y_{pre}+\Delta_{L,pre}$ and $\widehat{\Sigma}^{-1}=\Sigma_{pre}^{-1}+\Delta_{\Sigma,pre}^{-1}$,
the minimum of $||\widehat{L}\widehat{\Sigma}^{-1}-Y_{pre}\Sigma_{pre}^{-1}||$ is for $\widehat{L}=\widehat{L}_{UNALCE}$ and $\widehat{\Sigma}=\widehat{\Sigma}_{UNALCE}$.

%$$||(\widehat{S}^{-1/2}-Z_{pre}^{-1/2})(\widehat{S}^{-1/2}+Z_{pre}^{-1/2})[(I_p+\widehat{L}\widehat{S}^{-1})^{-1/2}-(I_p+Y_{pre}Z_{pre}^{-1/2})]$$
%$$[(I_p+\widehat{L}\widehat{S}^{-1})^{-1/2}+(I_p+Y_{pre}Z_{pre}^{-1/2})](\widehat{L}\widehat{S}^{-1}-Y_{pre}Z_{pre}^{-1})||.$$
%Since $(S^{1/2}+Z_{pre}^{1/2})=(S^{1/2}-Z_{pre}^{1/2}+2 Z_{pre}^{1/2})$,
%$||S^{1/2}-Z_{pre}^{1/2}||+2 ||Z_{pre}^{1/2}||$ and
%$((I_p+\widehat{L}\widehat{S}^{-1})^{-1/2}+(I_p+Y_{pre}Z_{pre}^{-1/2}))=((I_p+\widehat{L}\widehat{S}^{-1})^{-1/2}-(I_p+Y_{pre}Z_{pre}^{-1/2})+2 (I_p+Y_{pre}Z_{pre}^{-1/2}))$,
%$||(I_p+\widehat{L}\widehat{S}^{-1})^{-1/2}-(I_p+Y_{pre}Z_{pre}^{-1/2})||+2 ||(I_p+Y_{pre}Z_{pre}^{-1/2})||$
%, the thesis follows from the optimality properties of $\widehat{L}_{UNALCE}\widehat{S}_{UNALCE}^{-1}$,
%and the ones of $\widehat{S}_{UNALCE}^{-1/2}$ and $(I_p+\widehat{L}_{UNALCE}\widehat{S}_{UNALCE}^{-1})^{-1/2}$,
%which rely on the monotonicity of the variance of their estimated eigenvalues around their mean.

We can extend the validity of the proved optimality to $\widehat{f}_{k,T}-f_{k,T}=\widehat{B}'\widehat{\Sigma}^{-1}x_k-B'\Sigma^{*-1}x_k$
recalling that $||\widehat{f}_{k,T}-f_{k,T}|| \leq ||\widehat{B}'\widehat{\Sigma}^{-1}-B'\Sigma^{*-1}|| \times ||x_k||$ and
noting that
$$\widehat{B}\widehat{f}_{k,T}-Bf=\widehat{B}{\widehat{B}'\widehat{\Sigma}^{-1}}-BB'\Sigma^{*-1}=
\widehat{B}{\widehat{B}'\widehat{\Sigma}^{-1}}-BB'\Sigma^{*-1}+\widehat{B}B'\Sigma^{*-1}-\widehat{B}B'\Sigma^{*-1}=$$
$$=\widehat{B}(\widehat{B}'\widehat{\Sigma}^{-1}-B'\Sigma^{*-1})+(\widehat{B}-B)B'\Sigma^{*-1}.$$
%Since $\widehat{B}=O(p^{\alpha/2})$ and $\lambda_{i}(B)=\sqrt{\lambda_i(L^*)}$, $i=1,\ldots,r$, the eigenvalues of $\widehat{B}$ have the same optimality property for UNALCE, causing $(\widehat{B}-B)$ to be the minimum possible. Therefore, since $\widehat{B}{\widehat{B}'\widehat{\Sigma}^{-1}}-BB'\Sigma^{*-1}$ is the minimum possible,
%also $\widehat{B}'\widehat{\Sigma}^{-1}-B'\Sigma^{*-1}=\widehat{f}_{Thompson}-f$  must be optimal for UNALCE.
In fact, conditioning on ${B}_{pre}$ and $\Sigma_{pre}$, $||\widehat{B}'\widehat{\Sigma}^{-1}-B'\Sigma^{*-1}||$ must be minimum for $\widehat{B}=\widehat{B}_{UNALCE}$ and $\widehat{\Sigma}=\widehat{\Sigma}_{UNALCE}$ because $||\widehat{B}-B||$ and $||\widehat{B}\widehat{f}_{k,T}x_k^{-1}-Bfx_k^{-1}||$
are minimum for those matrices.
%for Corollary \ref{Coroll_B} and the previous thesis.

\section{Algebraic and parametric properties of matrix error estimates}\label{diag_trace}

We briefly examine the behaviour of $\widehat{L}_{POET}$ and $\widehat{S}_{POET}$
with respect to the projection operator $\mathbb{P}$ and the reference norm $g_\gamma$.
%(see \cite{fan2013large}).
Suppose that $\psi$ converges to $0$.
%$\frac{p^\alpha}{\sqrt{n}}$ diverges with $p=o(n^{2\alpha})$, which causes
In that case, the consistency of those estimates is also guaranteed,
i.e. $\widehat{L}_{POET}$ and $\widehat{S}_{POET}$ belong to $\mathcal{M}$.
However, we know from \cite{luo2011high} that $g_{\gamma}(\widehat{\Sigma}-\Sigma^{*},\widehat{\Sigma}-\Sigma^{*}))\leq 11 \psi$ would cause $g_\gamma(\mathcal{P}_{\mathcal{Y}^\perp}(\widehat{S}+\widehat{L}-\Sigma_n) < \psi$:
therefore, $\widehat{L}_{POET}$ and $\widehat{S}_{POET}$
would coincide with $\widehat{L}_{ALCE}$ and $\widehat{S}_{ALCE}$.
As a consequence, as $\psi$ is far from zero we have $||\mathbb{P}_{T'^\perp} (\widehat{L}_{POET}-L^*)||>\xi(T)\psi$
or $g_\gamma(\mathcal{P}_{\mathcal{Y}^\perp}(\widehat{S}_{POET}+\widehat{L}_{POET}-\Sigma_n)) > \psi$, which means $g_{\gamma}(\widehat{\Sigma}_{POET}-\Sigma^{*},\widehat{\Sigma}_{POET}-\Sigma^{*})> 11 \psi$.
As $\psi$ converges to $0$, instead, both POET and UNALCE estimates converge to the ALCE ones. %${L}_{ALCE}$ and $\widehat{S}_{ALCE}$

Concerning the trace of estimates, we note that the traces of $\widehat{\Sigma}_{UNALCE}$
and $\widehat{\Sigma}_{ALCE}$ differ from the trace of $\Sigma_n$, due to the use in the solution algorithm of the accelerated optimization scheme of \cite{nesterov2013gradient} (otherwise the equality would hold).
Since by definition $diag(\widehat{L}_{UNALCE})+diag(\widehat{S}_{UNALCE})=diag(\widehat{\Sigma}_{ALCE})$,
it follows instead that $trace(\widehat{\Sigma}_{ALCE})=trace(\widehat{\Sigma}_{UNALCE})$ and
$diag(\widehat{L}_{UNALCE}-\widehat{L}_{ALCE})=-diag(\widehat{S}_{UNALCE}-\widehat{S}_{ALCE})$.
This leads to the following equality $$||diag(\widehat{S}_{UNALCE}-\widehat{S}_{ALCE})||^2_{Fro}=||diag(\widehat{L}_{UNALCE}-\widehat{L}_{ALCE})||^2_{Fro}=
\sum_{i=1}^{p} (\widehat{{L}}_{UNALCE,ii}~-~\widehat{{L}}_{ALCE,ii})^2$$
and the following inequality
$$||diag(\widehat{L}_{UNALCE}-\widehat{L}_{ALCE})||^2_{Fro}\leq tr(\widehat{L}_{UNALCE}-\widehat{L}_{ALCE})^2=||\widehat{U}_{ALCE} \Lambda I_r \widehat{U}_{ALCE}'||^2_{Fro}=r\lambda^2.$$

As we know from \cite{farne2020large} that $||diag(\widehat{S}_{UNALCE}-\widehat{S}_{ALCE})||^2_{Fro}=||(\widehat{S}_{UNALCE}-\widehat{S}_{ALCE})||^2_{Fro}$,
if follows that
$$0<||(\widehat{S}_{UNALCE}-S^{*})||^2_{Fro}-||(\widehat{S}_{ALCE}-S^{*})||^2_{Fro}\leq r\lambda^2$$ and
$$0<||diag(\widehat{L}_{UNALCE}-L^{*})||^2_{Fro}-||diag(\widehat{L}_{ALCE}-L^{*})||^2_{Fro} \leq r\lambda^2,$$
which leads to the following corollary. %also holds
\begin{Coroll}
Conditionally on $Y_{pre}$ and $Z_{pre}$,
$$trace((\widehat{L}_{UNALCE}-L^*)^2)-trace((\widehat{L}_{ALCE}-L^*)^2)\leq trace(\widehat{L}_{UNALCE}-\widehat{L}_{ALCE})^2=r\psi^2,$$
%while $$trace((\widehat{\Sigma}_{UNALCE}-\Sigma^*)^2)-trace((\widehat{\Sigma}_{ALCE}-\Sigma^*)^2)=trace(\widehat{\Sigma}_{UNALCE}-\widehat{\Sigma}_{ALCE})^2=0$$
$$trace((\widehat{S}_{UNALCE}-S^*)^2)-trace((\widehat{S}_{ALCE}-S^*)^2)\leq trace(\widehat{S}_{ALCE}-\widehat{S}_{UNALCE})^2=r\psi^2.$$
\end{Coroll}
%because conditionally on $L_{pre}$ and $S_{pre}$

\bibliographystyle{chicago}%natbib
\bibliography{factors_bib}

\end{document}